%% file: main.tex
\documentclass[11pt,letterpaper,twoside,reqno,nosumlimits]{amsart}
\usepackage[numbers]{natbib}
\input{macro-3}
\allowdisplaybreaks
\begin{document}
\title[Profiles and global regularity for gCLM]
{Self-similar singularity and global regularity for the\\
hyper-convective generalized Constantin--Lax--Majda equation}
\author{Thomas Y. Hou}
\address{Applied and Computational Mathematics,
California Institute of Technology, Pasadena, CA, USA}
\email{hou@cms.caltech.edu}

\author{Xiang Qin}
\address{Applied and Computational Mathematics,
California Institute of Technology, Pasadena, CA, USA}
\email{xqin2@caltech.edu}

\author{Peicong Song}
\address{Applied and Computational Mathematics,
California Institute of Technology, Pasadena, CA, USA}
\email{psong2@caltech.edu}
\hypersetup{pdftitle={Self-similar singularity and global regularity for the hyper-convective generalized Constantin--Lax--Majda equation}}
\begin{abstract}
We study odd real-line solutions of the generalized
Constantin--Lax--Majda equation with transport parameter \(a\).
Near \(a=3/2\), we construct a local branch of compactly supported
rescaled profiles, smooth in the support interior, and prove
exponential stability in a weighted norm. The physical equation develops finite-time singularities for \(a<3/2\) and is global and decaying for \(a>3/2\) from smooth compactly supported initial data near the profiles, provided \(a\) is sufficiently close to \(3/2\).
We also prove global smoothness for odd compactly
supported smooth data that are nonpositive on the positive half-line and
satisfy a compressive condition at every interior zero, for
\(a\geq2\). Zero initial slope at the origin extends this
conclusion to \(a\geq1\). 
The smooth-data blowup result rules out a universal global-regularity
threshold below \(3/2\), and the local profile branch changes from
blowup to global decay at this value. The proofs combine dynamically
rescaled stability, Mellin-type estimates, and a Lagrangian convexity
mechanism.
\end{abstract}
\maketitle

\section{Introduction and main results}\label{sec:introduction}
We study the generalized Constantin--Lax--Majda equation (gCLM for short)
\begin{equation}\label{eq:model}
\omega_t+a u\omega_x=u_x\omega,\qquad
u_x=\bmH_{\bbR}\omega,\qquad u(0,t)=0,
\end{equation}
in which the real parameter \(a\) measures the strength of transport
relative to stretching. Physical space and time are denoted by \(x\)
and \(t\). We use the real-line Hilbert transform
\begin{equation}\label{eq:hilbert}
\bmH_{\bbR}g(x)=\frac1\pi\operatorname{p.v.}
\int_{\bbR}\frac{g(s)}{x-s}\idiff s.
\end{equation}
The case \(a=0\) is the Constantin--Lax--Majda model
\citep{ConstantinLaxMajda1985}. De Gregorio added the convective term
corresponding to \(a=1\) \citep{DeGregorio1990}, and Okamoto,
Sakajo, and Wunsch introduced the parameterized family
\citep{OkamotoSakajoWunsch2008}. We refer to the regime \(a>1\) as
hyper-convective.

{
The gCLM family is motivated by the vorticity formulation of the
three-dimensional incompressible Euler equations,
\begin{equation}\label{eq}
\partial_t\boldsymbol\omega
+\boldsymbol u\cdot\nabla\boldsymbol\omega
=\boldsymbol\omega\cdot\nabla\boldsymbol u,
\qquad
\boldsymbol\omega=\nabla\times\boldsymbol u,
\qquad
\nabla\cdot\boldsymbol u=0,
\end{equation}
where transport and
vortex stretching compete in the evolution of vorticity.
Equation~\eqref{eq:model} retains this competition in one dimension.
The term \(a u\omega_x\) transports vorticity, while \(u_x\omega\)
amplifies or attenuates it, and the order-zero relation
\(u_x=\bmH_{\bbR}\omega\) has the same scaling as the relation between
velocity gradient and vorticity in the Euler equations.  The
parameter \(a\) therefore provides a simple setting in which to study
how stronger transport affects singularity formation.

Finite-time singularity formation for the gCLM family is known in a
number of parameter and regularity regimes.  The CLM equation at
\(a=0\) admits an explicit finite-time singularity
\citep{ConstantinLaxMajda1985}, while Castro and C\'ordoba established
smooth-data blowup for negative transport parameters
\citep{CastroCordoba2010}.  On the real line, Elgindi and Jeong
constructed smooth self-similar blowup for sufficiently small
\(|a|\), and Elgindi, Ghoul, and Masmoudi subsequently proved
stability of these profiles under suitable weighted odd perturbations
\citep{ElgindiJeong2020Effects,ElgindiGhoulMasmoudi2021}.  Exact
collapsing solutions at \(a=1/2\) also arise in the work of Chen and
of Lushnikov, Silantyev, and Siegel
\citep{Chen2020,LushnikovSilantyevSiegel2021}.  At the De Gregorio
parameter \(a=1\), Chen, Hou, and Huang proved asymptotically
self-similar blowup from smooth compactly supported real-line data
using computer-assisted estimates \citep{ChenHouHuang2021}, while
Huang, Tong, and Wei constructed infinitely many exact compactly
supported self-similar profiles \citep{HuangTongWei2023}.  More
generally, interiorly smooth self-similar profiles are known for every
\(a\leq1\) \citep{HuangQinWangWei2024}.  At lower regularity,
self-similar singularity formation extends to a much wider range of
transport parameters when the H\"older exponent is sufficiently small
\citep{ElgindiJeong2020Effects,Zheng2023}.
Recent work of Xu analyzes the spectrum of the CLM collapse profile
and, under explicit profile hypotheses, the essential spectra of gCLM
linearizations for \(0<a<1\), emphasizing their dependence on the
operator domain \citep{Xu2026CLMSpectrum,Xu2026GCLMEssential}.

Closely related to these blowup results is the equally important
question of whether sufficiently strong transport suppresses
finite-time singularity formation, and how the answer depends on the
underlying domain.
 On the real line, the picture remains substantially incomplete.
  Lushnikov, Silantyev, and Siegel numerically found finite-time blowup up to
\(a=1\), with a change near \(a\approx0.6890665\) from focusing to
expanding self-similar behavior
\citep{LushnikovSilantyevSiegel2021}.  For stronger advection, they
reported numerical evidence of global evolution for the smooth
data tested at \(a\gtrsim1.3\), while their computations for
\(1<a\lesssim1.3\) remained inconclusive at large times.

On the circle, Jia, Stewart, and \v Sver\'ak established nonlinear stability near a
manifold of De Gregorio equilibria, while Lei, Liu, and Ren proved
weighted stability near the ground state \(-\sin\theta\)
\citep{JiaStewartSverak2019,LeiLiuRen2020}.  For odd data close to
this ground state, Chen showed finite-time blowup for \(a<1\)
sufficiently close to one and global decay for \(a>1\) sufficiently
close to one \citep{Chen2021}.  Further global regularity results for
periodic odd signed data impose quantitative conditions at the origin
and on the zero set \citep{Chen2025Regularity}, while related work
includes a BMO continuation criterion at \(a=2\) and steady or
traveling-wave solutions
\citep{Wunsch2011,OkamotoSakajoWunsch2014}.
More recently, Hou, Qin, and Song constructed periodic steady states
and rescaled profiles with a selected finite or infinite period, together
with sharp regularity and endpoint asymptotics \citep{Profiles}. They also established global
regularity under sign, origin, and interior-zero conditions, a blowup
criterion throughout \(0<a<1\), and realization of focusing or finite-period
profiles by smooth periodic blowup \citep{CompanionDynamics}.
These periodic results, together with the real-line numerical evidence,
raise the question of whether there is a critical transport strength
separating finite-time blowup from global evolution for smooth
real-line data.

The dynamics also depend strongly on the regularity and structure of
the data.  Huang, Qin, and Wang obtained multiscale blowup asymptotics
for the CLM equation for data with a degenerate zero at the origin
\citep{HuangQinWang2025}.
Huang, Tong, and Wang subsequently constructed explicit singular
profiles for \(a<0\), proved convergence to a singular outer profile
at \(a=0\) for suitable degenerate data, and established traveling
waves for \(a\leq1\), alongside numerical investigations of the
associated blowup scenarios \citep{HuangTongWang2026}.  Dissipative variants likewise reveal a
marked distinction between the line and the circle, with small
periodic data remaining globally regular in regimes where exact
real-line solutions can still develop finite-time singularities
\citep{Chen2020,AmbroseLushnikovSiegelSilantyev2024}.
Silantyev, Lushnikov, Siegel, and Ambrose subsequently obtained exact
periodic solutions with dissipation for \(a=0\) and \(a=1/2\), and
used them to compare self-similar blowup with small-data global
existence \citep{SilantyevLushnikovSiegelAmbrose2025}.
Related mechanisms also appear in other reduced fluid models.  The Hou--Luo
model admits a computer-assisted focusing-blowup construction from
smooth data \citep{ChenHouHuang2022HouLuo}, while exact
weak-advection Hou--Li profiles have recently been constructed under
periodic and real-line boundary conditions
\citep{HouQinWang2026HouLi}. Homogeneous SQG reductions lead to angular transport equations closely
related to the gCLM model at \(a=2\), with the same relative strength
of transport and stretching but a different nonlocal velocity law
\citep{ElgindiJeong2020Symmetries}.  These results show that the
transport parameter alone does not determine the dynamics and that
the domain, regularity, sign structure, and order of vanishing of the
data can all play essential roles.

Recent developments for the three-dimensional Euler equations further
emphasize this connection.  Elgindi \citep{Elgindi2021Euler} established
finite-time singularity formation on the whole space \(\bbR^3\) for
\(C^{1,\alpha}\) velocity fields with sufficiently small \(\alpha>0\).
Chen and Hou \citep{ChenHou2025Euler} established finite-time
singularity formation from smooth finite-energy initial data in an
axisymmetric setting with boundary.  In 2026, Shkoller
constructed finite-time blowup in the axisymmetric no-swirl class for
velocity in \(C^{1,\alpha}\cap L^2\) for every
\(0<\alpha<1/3\) \citep{Shkoller2026Euler}.  Chen constructed exact
self-similar profiles and asymptotically self-similar blowup in the
same regularity range, together with a one-dimensional limiting
profile as \(\alpha\) approaches \(1/3\)
\citep{Chen2026EulerI,Chen2026EulerII}, while Shao, Wei, Zhang, and
Zhang obtained related self-similar constructions in dimensions
\(d\geq3\) below the corresponding \(C^{1,1-2/d}\) threshold
\citep{ShaoWeiZhangZhang2026Euler}.  Smooth-forcing Euler blowup
constructions have also appeared \citep{AlpogeBuckmaster2026Euler},
and a recent construction gives finite-time blowup for the unforced
three-dimensional Euler equations from smooth compactly supported
divergence-free initial velocity \citep{OpenAI2026Euler}.  Parallel
progress has occurred for three-dimensional Navier--Stokes, where a
smooth-forcing construction produces finite-time singularity
formation on both \(\bbR^3\) and the periodic three-torus while the
kinetic energy remains bounded \citep{OpenAI2026NavierStokes}.  This
result concerns a prescribed smooth body force and is distinct from
the unforced Navier--Stokes Cauchy problem.  These higher-dimensional
developments reinforce the usefulness of reduced equations for
isolating singularity mechanisms, while also emphasizing the roles of
forcing, regularity, geometry, and data class.

}

In the present setting, we focus on smooth compactly supported
unforced solutions of the real-line gCLM equation in the
hyper-convective regime.
We seek finite-time singularities that are asymptotically self-similar
in the sense that, for suitable positive time-dependent scaling
factors \(C_\omega\) and \(C_l\),
\[
 \omega(x,t)
 =C_\omega
 \left(\bar\Omega(C_l x)+o(1)\right), \qquad C_\omega\longrightarrow+\infty
\text{ as }t\to T,
\]
in an appropriate topology.  The precise normalization of the
scaling factors and the convergence topology are specified below,
while \(\bar\Omega\) is the corresponding self-similar profile.
Our first result constructs a stable branch
of compactly supported profiles near \(a=3/2\).  For \(a<3/2\)
sufficiently close to \(3/2\), its stability basin contains smooth
compactly supported data whose solutions develop finite-time
asymptotically self-similar singularities.  Thus smooth-data blowup
persists strictly into the hyper-convective regime.  We also establish
global regularity for a class of odd one-sided signed data when
\(a\geq2\), subject to a condition at each interior zero, and extend
the result to \(a\geq1\) when the initial slope at the origin
vanishes.  The local branch and the signed-data theorem describe complementary
regimes: the former resolves the dynamics near a specified profile,
whereas the latter applies throughout the class defined below.

The profile construction uses a real rescaled coordinate \(y\in[-1,1]\).
Uppercase \(\Omega,U\) denotes the rescaled vorticity and velocity,
and \(\tau\) denotes rescaled time.
A bar identifies a stationary rescaled field, so
\(\bar\Omega_a,\bar U_a\) denote the profile family and its normalized
velocity. Stationarity here is in rescaled time; the corresponding
physical solution can grow, decay, or change its spatial scale.
The reference member is
\begin{equation}\label{eq:reference}
\bar\Omega_{3/2}(y)=-y\sqrt{1-y^2}\,\mathbf1_{\{|y|\leq1\}}.
\end{equation}
Here \(\mathbf1\) denotes the indicator of the indicated set.
For odd perturbations supported in this interval, the weighted norm
and differential operator are defined by
\begin{equation}\label{eq:real-spaces}
 \|g\|_{\bbH_y}^2=\frac1{2\pi}\int_{-1}^{1}
 \frac{\sqrt{1-y^2}}{y^2}|g_y|^2\idiff y,
 \qquad \calD_y g=y(1-y^2)g_y.
\end{equation}
The space \(\bbH_y\) is the completion of smooth interior odd
functions that vanish near the support endpoints and have zero slope
at the origin. Its graph space \(\bbX_y\) consists of the elements
with \(\calD_y g\in\bbH_y\), with norm
\(\|g\|_{\bbX_y}=\|g\|_{\bbH_y}+\|\calD_y g\|_{\bbH_y}\).
The trace, maximal-domain interpretation, and equivalent circle
realization are proved in Sec.~\ref{sec:coordinates}. 
We can now state the first main result, which describes the local
profile branch near the hyper-convective parameter \(a=3/2\) and its
consequences for the physical equation.

{
\begin{theorem}[Local profile branch and smooth-data dynamics]
\label{thm:profiles}
There are constants
\(\delta_{\rm br},r_{\rm br},\varepsilon_{\rm st},C_{\rm br}>0\),
independent of \(a\), such that for
\(|a-3/2|<\delta_{\rm br}\) the following statements hold.

\begin{enumerate}

\item[(i)] \emph{Profile branch.}
There is an odd normalized profile \(\bar\Omega_a\), supported in
\([-1,1]\), satisfying
\begin{equation}\label{eq:profile}
 (\bar c_{l,a}y+a \bar U_a)(\bar\Omega_a)_y
 =(\bar c_{\omega,a}+(\bar U_a)_y)\bar\Omega_a,\qquad
 (\bar U_a)_y=\bmH_{\bbR}\bar\Omega_a,\qquad \bar U_a(0)=0,
\end{equation}
where its normalization and rates obey
\begin{equation}\label{eq:profile-rates}
\begin{aligned}
 &(\bar\Omega_a)_y(0)=-1,\qquad
 \bar c_{l,a}=-a \bar U_a(1)<0,\\
 &\bar c_{\omega,a}=\bar c_{l,a}+(a-1)(\bar U_a)_y(0),\qquad
 \operatorname{sgn}\bar c_{\omega,a}=\operatorname{sgn}(a-3/2).
\end{aligned}
\end{equation}
It is the unique stationary profile in the affine class
\(\bar\Omega_{3/2}+\bbX_y\) within distance \(r_{\rm br}\) of
\(\bar\Omega_{3/2}\), and it satisfies the estimate
\[
 \|\bar\Omega_a-\bar\Omega_{3/2}\|_{\bbX_y}
 \leq C_{\rm br}|a-3/2|.
\]
The profile is smooth on \((-1,1)\), and its zero extension solves
\eqref{eq:profile} distributionally on \(\bbR\).

\item[(ii)] \emph{Corresponding physical solutions.}
Each profile \(\bar\Omega_a\) generates an exact solution of the
physical gCLM equation \eqref{eq:model}.  At \(a=3/2\), the profile equals \(\bar\Omega_{3/2}\), and its rates are
\(\bar c_{l,3/2}=-1/4\) and \(\bar c_{\omega,3/2}=0\). It generates the exact
solution \(\omega(x,t)=\bar\Omega_{3/2}(e^{-t/4}x)\).
For \(a\ne3/2\), the exact solution generated by \(\bar\Omega_a\) is
\begin{equation}\label{eq:exact-solution}
 \omega(x,t)=\frac1{1+\bar c_{\omega,a}t}
 \bar\Omega_a\bigl((1+\bar c_{\omega,a}t)^{\bar c_{l,a}/\bar c_{\omega,a}}x\bigr),
\end{equation}
on the interval specified by \(1+\bar c_{\omega,a}t>0\). Thus it develops a finite-time singularity for \(a<3/2\) and is global and decaying for
\(a>3/2\).

\item[(iii)] \emph{Smooth-data dynamics.}
For every such \(a\), the \(\bbH_y\)-neighborhood
\begin{equation}\label{eq:basin}
 \omega_0-\bar\Omega_a\in\bbX_y,\qquad
 \|\omega_0-\bar\Omega_a\|_{\bbH_y}<\varepsilon_{\rm st},
\end{equation}
contains odd initial data
\(\omega_0\in C_c^\infty(\bbR)\) supported strictly inside
\((-1,1)\).  For every such smooth initial datum, the solution of the
physical equation \eqref{eq:model} converges after rescaling to
\(\bar\Omega_a\) as in \eqref{eq:stability}.
If \(a<3/2\), it develops a finite-time asymptotically self-similar
singularity with unbounded vorticity amplitude.
If \(a>3/2\), it is global and smooth, and its amplitude satisfies
\[
 \lim_{t\to\infty}t\,\|\omega(\cdot,t)\|_{L^\infty(\bbR)}
 =\frac{\|\bar\Omega_a\|_{L^\infty(\bbR)}}{\bar c_{\omega,a}}.
\]

\end{enumerate}
\end{theorem}
}

The blowup solutions are spatially expanding: since
\(\bar c_{l,a}<0\), their characteristic length scale \(C_l^{-1}\)
tends to infinity as the vorticity amplitude diverges; see
\eqref{eq:scale-asymptotics}.

We next turn from dynamics near the profile branch to
global regularity in the hyper-convective regime.  The result applies
to the following class of odd signed data.
\begin{definition}[Admissible signed data]\label{def:data}
For a nonzero odd function \(g\in C_c^\infty(\bbR)\), let
\(R_g=\sup\{x>0:g(x)\ne0\}\). Define the signed-data class by
\begin{equation}\label{eq:compressive-zero}
 \bbG=\{0\}\cup
 \left\{g\in C_c^\infty(\bbR)\setminus\{0\}:
 \begin{array}{l}
 g\text{ is odd and }g(x)\leq0\text{ for }x>0,\\
 g(z)=0\ \Longrightarrow\ (\bmH_{\bbR}g)(z)\leq0
       \text{ for every }z\in(0,R_g)
 \end{array}\right\}.
\end{equation}
We call the last requirement the \emph{compressive-zero condition}.
It is vacuous when \(g<0\) throughout \((0,R_g)\).
\end{definition}

\begin{theorem}[Signed-data global regularity]\label{thm:global}
For a parameter \(a\) in the range \(a\geq2\) and initial data
\(\omega_0\) in \(\bbG\), the solution of
\eqref{eq:model} exists for all \(t\geq0\) and remains smooth.
The same conclusion holds in the range \(a\geq1\) under the
additional condition \(\omega_{0x}(0)=0\).
\end{theorem}

Several parts of the argument are motivated by earlier work.  The
circle formulation and weighted framework are inspired in part by the
De Gregorio analysis of \citep{LeiLiuRen2020}, while the rescaling and
perturbative stability arguments are related to those developed in
\citep{ChenHouHuang2021,Chen2021}.  The Lagrangian viewpoint in the
global regularity argument is also connected to the geometric and
solar formulations in \citep{BauerKolevPreston2016,Preston2023}.
The weighted spaces, profile construction, Mellin estimates, positive
force identity, and noncollapse argument used here are developed
directly for the present problem.

The remainder of the paper is organized as follows.
Section~\ref{sec:coordinates} develops the correspondence between the
compactly supported real-line profile and its circle representation,
where the reference profile becomes \(-\sin\theta\).
Section~\ref{sec:linear} establishes the linear theory on the weighted
graph space and constructs the local profile branch by the implicit
function theorem, including the sign change of the amplitude scaling
rate.  Section~\ref{sec:stability} proves nonlinear stability of this
branch and reconstructs the corresponding physical solutions,
including smooth compactly supported data in the blowup and decay regimes.
Section~\ref{sec:global} proves the global regularity results in the
hyper-convective regime by combining moment estimates with a
Lagrangian noncollapse argument.  Appendix~\ref{app:rows} contains
the analytic matrix estimates used in the profile construction.

\section{Coordinates and weighted domains}\label{sec:coordinates}
In this section, we set up the rescaled and functional framework used in
both the profile construction and the nonlinear stability analysis.
Since the reference profile at \(a=3/2\) is compactly supported, we
seek the nearby profile branch within the same compactly supported
class.  We first introduce the normalized rescaling, fixing the origin
slope and the support endpoints at \(y=\pm1\).  With the support thus
normalized to \([-1,1]\), an exact change of variables transforms the
real-line problem to the circle, where the reference profile becomes
\(-\sin\theta\) and the Hilbert transform admits a convenient
conjugacy.  Finally, we introduce the weighted Hilbert and graph spaces
adapted to the origin and support endpoints, and identify their traces,
basis, and maximal domains.  These constructions provide the common
setting for the stationary and evolutionary arguments below.

The real rescaled fields \(\Omega,U\) are related to the physical
fields by
\begin{equation}\label{eq:rescaling}
 \Omega(y,\tau)=C_\omega(\tau)^{-1}
 \omega(C_l(\tau)^{-1}y,t(\tau)),\qquad
 U(y,\tau)=C_\omega(\tau)^{-1}C_l(\tau)
 u(C_l(\tau)^{-1}y,t(\tau)),
\end{equation}
where the positive scaling functions are defined by
\begin{equation}\label{eq:scales}
 C_\omega(\tau)=\exp\left(-\int_0^\tau c_\omega(s)\idiff s\right),
 \quad C_l(\tau)=\exp\left(\int_0^\tau c_l(s)\idiff s\right),
 \quad t(\tau)=\int_0^\tau C_\omega(s)^{-1}\idiff s.
\end{equation}
Differentiating \eqref{eq:rescaling}, using
\(C_\omega'/C_\omega=-c_\omega\), \(C_l'/C_l=c_l\), and
\(t'=C_\omega^{-1}\), gives the equation
\begin{equation}\label{eq:rescaled}
 \Omega_\tau+(c_l y+aU)\Omega_y=(c_\omega+U_y)\Omega,
 \qquad U_y=\bmH_{\bbR}\Omega,\quad U(0,\tau)=0.
\end{equation}
We choose the normalization
\begin{equation}\label{eq:normalization}
 c_l=-aU(1),\qquad c_\omega=c_l+(a-1)U_y(0).
\end{equation}
The effective velocity vanishes at \(y=0,\pm1\). Differentiating
\eqref{eq:rescaled} at the origin gives the slope equation
\begin{equation}\label{eq:normalized-slope}
 \partial_\tau\Omega_y(0,\tau)
 =\bigl(c_\omega-c_l+(1-a)U_y(0,\tau)\bigr)\Omega_y(0,\tau)=0,
\end{equation}
where the last equality uses \eqref{eq:normalization}. Thus the
initial slope \(-1\) is preserved. The weighted spaces below retain
this normalization in the trace sense.
We now pass to circle coordinates, where the fixed support interval
is absorbed into the geometry and the reference profile takes the
simpler form \(-\sin\theta\).

The circle map and its action on a test function are given by
\begin{equation}\label{eq:circle-map}
 y=\sin(\theta/2),\qquad
 \widetilde\Omega(\theta,\tau)=2\Omega(\sin(\theta/2),\tau),
 \qquad (\calT g)(\theta)=2g(\sin(\theta/2)).
\end{equation}
In particular, the identity
\(\widetilde{\bar\Omega}_{3/2}=\calT\bar\Omega_{3/2}=-\sin\theta\)
defines the circle reference profile.
We write \(\bbT=\bbR/(2\pi\bbZ)\) and use the periodic Hilbert transform
\begin{equation}\label{eq:periodic-hilbert}
\bmH_{\bbT}g(\theta)
=\frac1{2\pi}\operatorname{p.v.}
\int_{-\pi}^{\pi}g(s)\cot\frac{\theta-s}{2}\idiff s.
\end{equation}

\begin{lemma}[Odd Hilbert conjugacy]\label{lem:conjugacy}
For an odd input \(g\) in \(L^2(\bbT)\), define
\(\calT^{-1}g\) on \((-1,1)\) by the formula
\(\calT^{-1}g(y)=g(2\arcsin y)/2\), and extend it by zero.
Then the identity
\begin{equation}\label{eq:conjugacy}
 (\bmH_{\bbR}\calT^{-1}g)(\sin(\theta/2))
 =\frac12(\bmH_{\bbT}g)(\theta)
\end{equation}
holds in \(L^2((-1,1),\idiff y)\) after the stated change of variables.
\end{lemma}
\begin{proof}
For smooth odd inputs, pair the positive and negative halves in both
principal values. With \(y=\sin(\theta/2)\) and \(r=\sin(s/2)\),
the elementary identity
\[
 \cot\frac{\theta-s}{2}-\cot\frac{\theta+s}{2}
 =\frac{2r\cos(s/2)}{y^2-r^2}
\]
and the relation \(\idiff r=\cos(s/2)\idiff s/2\) give
\[
 \begin{aligned}
 (\bmH_{\bbT}g)(\theta)
 &=\frac1\pi\operatorname{p.v.}\int_0^\pi
      \frac{r\cos(s/2)g(s)}{y^2-r^2}\idiff s\\
 &=\frac2\pi\operatorname{p.v.}\int_0^1
      \frac{r g(2\arcsin r)}{y^2-r^2}\idiff r
 =2(\bmH_{\bbR}\calT^{-1}g)(y).
 \end{aligned}
\]
This proves \eqref{eq:conjugacy} for smooth inputs. Approximate an odd \(L^2\)
input by smooth odd functions. Pullback to the interval is bounded
by the Jacobian identity
\(\idiff y=\cos(\theta/2)\idiff\theta/2\); the
\(L^2\) bounds for both Hilbert transforms therefore pass the
identity to the limit. Oddness is essential in the paired kernels.
\end{proof}

The conjugacy allows us to rewrite the normalized dynamics entirely on the circle. We next introduce the functionals that encode the normalization and effective transport.
For odd \(H^1(\bbT)\) inputs, the velocity trace and effective transport
operator are defined by
\begin{equation}\label{eq:circle-functionals}
 \mathfrak u_1(g)=\frac14\int_0^\pi\cos\frac s2
 (\bmH_{\bbT}g)(s)\idiff s,\qquad
 \mathfrak c_a(g)=-a\mathfrak u_1(g)+\frac{a-1}{2}(\bmH_{\bbT}g)(0),
\end{equation}
and
\begin{equation}\label{eq:P}
 (\calP g)(\theta)=\frac1{\cos(\theta/2)}
 \int_0^\theta\cos\frac s2
 \left(\frac12(\bmH_{\bbT}g)(s)-\mathfrak u_1(g)\right)\idiff s.
\end{equation}
The endpoint numerator at \(\pi\) vanishes by the definition of
\(\mathfrak u_1\). Direct differentiation gives
\begin{equation}\label{eq:P-ode}
 (\calP g)_\theta-\frac12\tan(\theta/2)\calP g
 =\frac12\bmH_{\bbT}g-\mathfrak u_1(g).
\end{equation}
For \(g=\widetilde\Omega\), this operator equals
\(2(U(y)-U(1)y)/\sqrt{1-y^2}\). Consequently the circle equation is
\begin{equation}\label{eq:circle}
 \widetilde\Omega_\tau=\calG(a,\widetilde\Omega),\qquad
 \calG(a,g)=-a(\calP g)g_\theta+
 \left(\mathfrak c_a(g)+\frac12\bmH_{\bbT}g\right)g.
\end{equation}
The reference values
\begin{equation}\label{eq:circle-reference}
 \bmH_{\bbT}\widetilde{\bar\Omega}_{3/2}=\cos\theta,
 \qquad\calP\widetilde{\bar\Omega}_{3/2}=\frac13\sin\theta,
 \qquad\mathfrak u_1(\widetilde{\bar\Omega}_{3/2})=\frac16
\end{equation}
give the stationary identity
\(\calG(3/2,\widetilde{\bar\Omega}_{3/2})=0\).
They also give \((\bar U_{3/2})_y=1/2-y^2\) and
\(\bar U_{3/2}=y/2-y^3/3\) on \(|y|\leq1\), proving the
reference normalization. This polynomial describes the velocity
on the support interval; the velocity is not extended by zero.

The weighted domain has two roles. Its origin weight fixes the
perturbation's zero slope, while its graph derivative follows a
transport field that vanishes at the support endpoints. The next
definition and lemma identify both traces and a common approximation
core, which will be needed for the closed linear operator and for
the nonlinear energy identities.

\begin{definition}[Circle spaces]\label{def:spaces}
Let \(\bbV\) be the real smooth odd periodic functions with
\(g_\theta(0)=0\). The Hilbert space \(\bbH_\theta\) is its
completion in the norm
\begin{equation}\label{eq:circle-norm}
 \|g\|_{\bbH_\theta}^2=\frac1{4\pi}\int_{-\pi}^{\pi}
 \frac{|g_\theta|^2}{\sin^2(\theta/2)}\idiff\theta.
\end{equation}
Its inner product \(\langle\cdot,\cdot\rangle_{\bbH_\theta}\)
is obtained by polarization.
The distributional operator and its graph norms are defined by
\begin{equation}\label{eq:graph}
 \calD_\theta g=\sin\theta\,g_\theta,\qquad
 \|g\|_{\bbX_\theta}=\|g\|_{\bbH_\theta}
 +\|\calD_\theta g\|_{\bbH_\theta},\qquad
 \|g\|_{\bbX_{\theta,2}}=\|g\|_{\bbX_\theta}
 +\|\calD_\theta^2g\|_{\bbH_\theta},
\end{equation}
on their maximal domains in \(\bbH_\theta\).
\end{definition}

\begin{lemma}[Traces, basis and maximal domain]\label{lem:domain}
The space \(\bbH_\theta\) embeds continuously in odd periodic
\(H^1\), with \(g(0)=g(\pi)=0\) and
\(g(\theta)/\theta\to0\) at the origin. For positive integer
indices \(i,j,k\), the functions
\begin{equation}\label{eq:basis}
 e_j(\theta)=\frac{\sin((j+1)\theta)}{j+1}
 -\frac{\sin(j\theta)}j,\qquad e_0=0,
\end{equation}
form an orthonormal basis, with the expansion and Parseval identity
\begin{equation}\label{eq:basis-expansion}
 g=\sum_{j\geq1}b_j e_j,\qquad
 \|g\|_{\bbH_\theta}^2=\sum_{j\geq1}|b_j|^2.
\end{equation}
The operator \(\calD_\theta\), with maximal domain
\(\bbX_\theta\subset\bbH_\theta\), is closed.  The finite span of
the functions \(e_j\) is a core for both graph domains in
\eqref{eq:graph}.  The map \(\calT\) is an isometry from
\(\bbH_y\) to \(\bbH_\theta\) and from \(\bbX_y\) to
\(\bbX_\theta\), and it intertwines their differential operators.
\end{lemma}
\begin{proof}
The weighted derivative \(g_\theta/\sin(\theta/2)\) identifies the
half-circle space with the closed subspace of \(L^2(0,\pi)\) satisfying
\begin{equation}\label{eq:compatibility}
 \int_0^\pi g_\theta\idiff\theta=0.
\end{equation}
Conversely, integration against \(\sin(\theta/2)\) reconstructs
\(g\) from any compatible weighted derivative. Cauchy--Schwarz gives
\(g(\theta)=o(\theta^{3/2})\) at zero, and ordinary \(H^1\)
control supplies both endpoint traces. The identity
\((e_j)_\theta/\sin(\theta/2)=-2\sin((j+1/2)\theta)\)
identifies the basis with the complete half-integer sine system;
compatibility removes its first mode \(\sin(\theta/2)\).

We next turn from the Hilbert-space structure to the maximal graph
domains and their approximation cores.  We first verify the closedness
of the maximal graph operator.  Suppose that \(g_n\in\bbX_\theta\)
converges to \(g\in\bbH_\theta\) and that
\(\calD_\theta g_n\) also converges in \(\bbH_\theta\).
The continuous embedding
\(\bbH_\theta\hookrightarrow H^1(\bbT)\) permits passage to the
distributional limit in the relation
\(\calD_\theta g_n=\sin\theta\,\partial_\theta g_n\).
The limit of \(\calD_\theta g_n\) is therefore
\(\sin\theta\,\partial_\theta g\) in the distributional sense.
Since this limit belongs to \(\bbH_\theta\), the definition of the
maximal domain places \(g\) in \(\bbX_\theta\), and the convergence
is \(\calD_\theta g_n\to\calD_\theta g\).  This proves the closedness of
\(\calD_\theta\) on \(\bbX_\theta\).

To identify a core,
we use the coordinate on \(0<\theta<\pi\),
\begin{equation}\label{eq:graph-coordinate}
 \zeta=\log\tan(\theta/2),\qquad
 \calD_\theta=\partial_\zeta,\qquad
 \idiff\theta=\operatorname{sech}\zeta\idiff\zeta,
\end{equation}
and the exact identity
\begin{equation}\label{eq:weighted-graph}
 \frac{(\calD_\theta g)_\theta}{\sin(\theta/2)}
 =\partial_\zeta\left(\frac{g_\theta}{\sin(\theta/2)}\right)
 +\frac{1+3\cos\theta}{2}\cdot\frac{g_\theta}{\sin(\theta/2)}.
\end{equation}
The last coefficient and its derivative are bounded. Conversely,
the weighted \(H^1\) condition on \(g_\theta/\sin(\theta/2)\)
places
\(\sqrt{\operatorname{sech}\zeta}\,g_\theta/\sin(\theta/2)\)
in ordinary \(H^1(\bbR)\), since the logarithmic derivative of
the weight is bounded. This function tends to zero at both infinities.
The identity
\[
 \calD_\theta g
 =\sqrt{\operatorname{sech}\zeta}\sin(\theta/2)
   \left(\sqrt{\operatorname{sech}\zeta}\,
                       \frac{g_\theta}{\sin(\theta/2)}\right)
\]
therefore gives zero endpoint traces for \(\calD_\theta g\).
Its weighted derivative satisfies the compatibility condition by
integration. Iteration proves the same assertion at the second
graph level. Thus the two maximal graph domains correspond exactly
to compatible weighted \(H^1\) and
\(H^2\) spaces with measure \(\operatorname{sech}\zeta\idiff\zeta\).
Cutting off and mollifying the weighted derivative in \(\zeta\)
converges in these spaces. Subtract its compatibility error times a
fixed compactly supported smooth function whose compatibility integral
is one. The error tends to zero, so reconstruction gives smooth
functions vanishing near \(0,\pi\) that converge in the graph norm.
Their sine partial sums converge in every \(C^\ell\) norm, where
\(\ell\) denotes a nonnegative integer differentiation order. Subtracting
the origin slope times \(\sin\theta\) preserves that convergence
and restores zero slope. These corrected polynomials are finite
combinations of the \(e_j\). 
For odd zero-slope functions,
\(C^3\) convergence controls the first graph norm, and \(C^4\)
convergence controls the second. This proves the core assertion
without assuming convergence of arbitrary sharp basis truncations.

Finally, the change of variables \(y=\sin(\theta/2)\) in
\eqref{eq:circle-norm} gives the norm identity
\[
 \|\calT g\|_{\bbH_\theta}=\|g\|_{\bbH_y}.
\]
Differentiation of the circle map gives the intertwining identity
\[
 \calD_\theta\calT g=\calT\calD_y g.
\]
The same two identities give the isometry of the corresponding graph
norms.  Density of the endpoint-cutoff cores on the two sides then
extends these identifications to the maximal domains.
\end{proof}

\section{The linear isomorphism and the profile branch}\label{sec:linear}

This section establishes the local stationary profile branch near
\(a=3/2\).  We begin with the linearized operator at the reference
profile \(-\sin\theta\), proving coercivity and invertibility on the
maximal graph domain.  The argument exploits a decomposition into a
tridiagonal De Gregorio part and a compact remainder, with the
required matrix estimates deferred to App.~\ref{app:rows}.  We then
establish the weighted product estimates and analyticity needed for
the nonlinear stationary problem.  These results allow us to apply
the implicit function theorem to construct the nearby profile branch.
A final transversality estimate determines the sign of the amplitude
scaling rate along the branch.

Linearizing \eqref{eq:circle} at \((3/2,-\sin\theta)\) gives
\begin{equation}\label{eq:linear}
 \calL g=-\frac12\sin\theta\,g_\theta+\frac12\cos\theta\,g
 +\frac32(\calP g)\cos\theta
 -\left(\mathfrak c_{3/2}(g)+\frac12\bmH_{\bbT}g\right)\sin\theta.
\end{equation}
Define the auxiliary operators by
\begin{align}
 \calR g&=3\calP g-\int_0^\theta(\bmH_{\bbT}g)(s)\idiff s,
 \label{eq:R}\\
 \calL_{\rm DG}g&=-\calD_\theta g+\cos\theta\,g
 -(\bmH_{\bbT}g)\sin\theta
 +\cos\theta\int_0^\theta(\bmH_{\bbT}g)(s)\idiff s,
 \label{eq:DG}\\
 \calK g&=\frac12\left((\calR g)\cos\theta
 -(\calR g)_\theta(0)\sin\theta\right).
 \label{eq:K}
\end{align}
The slope identity from \eqref{eq:P-ode} then proves the decomposition
\[
 \calL=\frac12\calL_{\rm DG}+\calK.
\]
Here \(\calL_{\rm DG}\) is the De Gregorio linearization. All its
needed properties are proved here.
To quantify the contribution of the compact remainder \(\calK\) to the
energy estimate, we write its matrix in the weighted basis.
Define the matrix entries, symmetric part, and absolute symmetric
row quantities by
\begin{align}
 K_{ij}&=\langle\calK e_j,e_i\rangle_{\bbH_\theta},\qquad
 K_{ij}^{\rm sym}=(K_{ij}+K_{ji})/2,\label{eq:matrix}\\
 \Gamma_i&=K_{ii}+\sum_{j\ne i}|K_{ij}^{\rm sym}|.\label{eq:rows}
\end{align}
The formulas in App.~\ref{app:rows} in particular prove convergence
of these row sums. 
The following proposition provides the main linear input for both the profile construction and the nonlinear stability analysis.

\begin{proposition}[Coercivity and invertibility]\label{prop:linear}
The operator \(\calL\), with domain \(\bbX_\theta\) in
\(\bbH_\theta\), is closed and onto, and it satisfies
\begin{equation}\label{eq:coercivity}
 \langle\calL g,g\rangle_{\bbH_\theta}
 \leq-\frac3{16}\|g\|_{\bbH_\theta}^2.
\end{equation}
Its inverse has the bounds
\begin{equation}\label{eq:inverse}
 \|\calL^{-1}g\|_{\bbH_\theta}\leq\frac{16}3\|g\|_{\bbH_\theta},
 \qquad
 \|\calL^{-1}g\|_{\bbX_\theta}
 \leq\left(\frac{44}3+\frac{32}3
 \|\calK\|_{\bbH_\theta\to\bbH_\theta}\right)\|g\|_{\bbH_\theta}.
\end{equation}
\end{proposition}
\begin{proof}
Product-to-sum identities and \(\bmH_{\bbT}\sin(k\theta)=-\cos(k\theta)\)
give
\[
 \calL_{\rm DG}\sin(k\theta)
 =-\frac{(k-1)^2}{2k}\sin((k+1)\theta)
 +\frac{k^2-1}{2k}\sin((k-1)\theta).
\]
Applying this identity to \(e_j\) yields
\begin{equation}\label{eq:DG-basis}
 \calL_{\rm DG}e_j=-\delta_{j+1}^{\rm DG}e_{j+1}
 -(\delta_{j+1}^{\rm DG}-\delta_j^{\rm DG})e_j
 +\delta_j^{\rm DG}e_{j-1},\qquad
 \delta_j^{\rm DG}=\frac{(j-1)^2(j+1)}{2j^2}.
\end{equation}

Recall that we write the basis
expansion \(g=\sum_{j\geq1}b_j e_j\). The off-diagonal terms cancel in the quadratic form of
\(\calL_{\rm DG}\), so the coefficients
\begin{equation}\label{eq:DG-energy}
 \mu_j^{\rm DG}
 =\frac{\delta_{j+1}^{\rm DG}-\delta_j^{\rm DG}}2
 =\frac14+\frac{j^2-j-1}{4j^2(j+1)^2}
\end{equation}
give
\[
 \left\langle\frac12\calL_{\rm DG}g,g\right\rangle_{\bbH_\theta}
 =-\sum_j\mu_j^{\rm DG}|b_j|^2
\]
on the finite core.  The first coefficient satisfies
\(\mu_1^{\rm DG}=3/16\), while
\(\mu_j^{\rm DG}\geq1/4\) for every \(j\geq2\).

We next control the contribution of the compact remainder
\(\calK\).  By the definition of the symmetric matrix entries,
we have
\[
 \langle\calK g,g\rangle_{\bbH_\theta}
 =\sum_i K_{ii}|b_i|^2
   +2\sum_{i<j}K_{ij}^{\rm sym}b_i b_j.
\]
Using \(2|b_i b_j|\leq |b_i|^2+|b_j|^2\) and the definition
\[
 \Gamma_i=K_{ii}+\sum_{j\ne i}|K_{ij}^{\rm sym}|,
\]
we obtain
\[
 \langle\calK g,g\rangle_{\bbH_\theta}
 \leq\sum_i\Gamma_i|b_i|^2.
\]
Lemma~\ref{lem:rows} proves that
\(\Gamma_1\leq0\) and \(\Gamma_j\leq1/16\) for every \(j\geq2\).
Consequently,
\[
 \langle\calL g,g\rangle_{\bbH_\theta}
 \leq-\frac3{16}|b_1|^2
      -\sum_{j\geq2}\left(\frac14-\frac1{16}\right)|b_j|^2
 \leq-\frac3{16}\|g\|_{\bbH_\theta}^2,
\]
which proves \eqref{eq:coercivity} on the finite core.

We next turn from the energy estimate to the closedness and surjectivity of the De Gregorio
part.  For this purpose, compute
\begin{equation}\label{eq:D-basis}
 \calD_\theta e_j=-\frac{j-1}2e_{j-1}
 +\frac12e_j+\frac{j+2}2e_{j+1}.
\end{equation}
The bounded remainder \(\calB_0=\calL_{\rm DG}+\calD_\theta\)
therefore has the matrix
\begin{equation}\label{eq:B0}
 \calB_0e_j=-\frac{j-1}{2j^2}e_{j-1}
 -\frac{j^2-j-1}{2j^2(j+1)^2}e_j
 +\frac{(j+2)(2j+1)}{2(j+1)^2}e_{j+1}.
\end{equation}
The three band operators satisfy the coefficient bounds
\[
 \sup_{j\geq1}\frac{j-1}{2j^2}\leq\frac18,\qquad
 \sup_{j\geq1}
 \left|\frac{j^2-j-1}{2j^2(j+1)^2}\right|\leq\frac18
,\qquad
 \frac{(j+2)(2j+1)}{2(j+1)^2}
 =1+\frac{j}{2(j+1)^2}\leq\frac98,
\]
which implies the bound
\[\|\calB_0\|_{\bbH_\theta\to\bbH_\theta}\leq \frac{11}{8}.\]
 Lemma~\ref{lem:domain} shows that \(\calD_\theta\) is closed on its
maximal graph domain \(\bbX_\theta\).  Since \(\calB_0\) is bounded
on \(\bbH_\theta\), the decomposition
\(\calL_{\rm DG}=-\calD_\theta+\calB_0\) shows that
\(\calL_{\rm DG}\), with domain \(\bbX_\theta\), is also closed. The preceding De Gregorio energy estimate, together with
Cauchy--Schwarz, gives
\[
 \frac38\|g\|_{\bbH_\theta}^2
 \leq
 -\langle\calL_{\rm DG}g,g\rangle_{\bbH_\theta}
 \leq
 \|\calL_{\rm DG}g\|_{\bbH_\theta}\|g\|_{\bbH_\theta},
\]
and hence the bound
\[\|g\|_{\bbH_\theta}\leq\frac83\|\calL_{\rm DG}g\|_{\bbH_\theta}.\]
If the images form a Cauchy sequence, this bound and
\(\calD_\theta g=\calB_0g-\calL_{\rm DG}g\) make the
preimages Cauchy in the graph norm. Thus its range is closed.

Next we prove surjectivity of \(\calL_{\rm DG}\).  Since its
range is already known to be closed, it suffices to prove that the
range is dense in \(\bbH_\theta\).  We therefore consider a function \(g\in\bbH_\theta\) orthogonal to the range. Orthogonality to the image of the first basis vector, together with
\eqref{eq:DG-basis}, gives the relation \(b_2=-b_1\).
More generally,
orthogonality to the image of each remaining basis vector gives the
recurrence
\[
 \delta_{j+1}^{\rm DG}(b_{j+1}+b_j)
 =\delta_j^{\rm DG}(b_j+b_{j-1}),
 \qquad j\geq2.
\]
Induction therefore gives
\(b_{j+1}=-b_j\).
Since the coefficient sequence
\((b_j)_{j\geq1}\) is square summable, this relation forces
\(b_j=0\) for every \(j\).  Thus the orthogonal complement of
\(\operatorname{Ran}\calL_{\rm DG}\) is trivial, and the range is
dense.
The range is both dense and closed, and hence equals
\(\bbH_\theta\).  Together with the preceding coercive estimate, this
proves that \(\calL_{\rm DG}\) is an isomorphism from
\(\bbX_\theta\) onto \(\bbH_\theta\).  The decomposition
\(\calL_{\rm DG}=-\calD_\theta+\calB_0\), together with
\(\|\calB_0\|_{\bbH_\theta\to\bbH_\theta}\leq11/8\) and the
preceding low-norm estimate, gives the graph-derivative bound
\[
\begin{aligned}
\|\calD_\theta g\|_{\bbH_\theta}
&\leq
\|\calB_0g\|_{\bbH_\theta}
+
\|\calL_{\rm DG}g\|_{\bbH_\theta}
\leq
\frac{11}{8}\|g\|_{\bbH_\theta}
+
\|\calL_{\rm DG}g\|_{\bbH_\theta}\\
&\leq
\left(
\frac{11}{8}\cdot\frac83+1
\right)
\|\calL_{\rm DG}g\|_{\bbH_\theta}
=
\frac{14}{3}
\|\calL_{\rm DG}g\|_{\bbH_\theta}.
\end{aligned}
\]
Combining this bound with
\(\|g\|_{\bbH_\theta}\leq(8/3)\|\calL_{\rm DG}g\|_{\bbH_\theta}\)
and the definition of the graph norm gives the estimate
\[
\|g\|_{\bbX_\theta}=\|g\|_{\bbH_\theta}+\|\calD_\theta g\|_{\bbH_\theta}
\leq
\left(
\frac83+\frac{14}{3}
\right)
\|\calL_{\rm DG}g\|_{\bbH_\theta}
=
\frac{22}{3}
\|\calL_{\rm DG}g\|_{\bbH_\theta},
\]
and hence the inverse estimate
\[
 \|\calL_{\rm DG}^{-1}g\|_{\bbX_\theta}
 \leq\frac{22}{3}\|g\|_{\bbH_\theta}.
\]

It remains to pass from the De Gregorio part to the full linearization.
The Hilbert--Schmidt estimate \eqref{eq:K-HS} in
Appendix~\ref{app:rows} makes \(\calK\) compact on
\(\bbH_\theta\). Since \(\calL_{\rm DG}^{-1}\) is bounded, the operator
\[
 \calL\calL_{\rm DG}^{-1}
 =\frac12\Id+\calK\calL_{\rm DG}^{-1}
\]
is a compact perturbation of an invertible operator and hence is
Fredholm of index zero.  In particular, its range is closed and its
kernel and cokernel have the same finite dimension.
The decomposition
\[
 \calL=-\frac12\calD_\theta+\frac12\calB_0+\calK
\]
gives the closed realization of \(\calL\) on the domain
\(\bbX_\theta\), since \(\calD_\theta\) is closed there and
\(\calB_0/2+\calK\) is bounded on \(\bbH_\theta\).  The same finite
span of the basis functions remains a core under this bounded
perturbation.  Hence the coercive estimate
\eqref{eq:coercivity}, proved above on the finite core, extends to
every \(g\in\bbX_\theta\).
The extended coercive estimate implies triviality of the kernel.
Indeed, for \(g\in\ker\calL\), the estimate
\[
 0
 =\langle\calL g,g\rangle_{\bbH_\theta}
 \leq-\frac3{16}\|g\|_{\bbH_\theta}^2
\]
forces \(g=0\).  Since
\(\calL\calL_{\rm DG}^{-1}\) is Fredholm of index zero and
\(\calL_{\rm DG}^{-1}\) is an isomorphism, \(\calL\) has index zero
as well.  The trivial kernel therefore implies a trivial cokernel,
and hence
\(\operatorname{Ran}\calL=\bbH_\theta\).

The same coercive estimate, together with Cauchy--Schwarz, gives
the low-norm bound
\[
 \frac3{16}\|g\|_{\bbH_\theta}^2
 \leq-\langle\calL g,g\rangle_{\bbH_\theta}
 \leq
 \|\calL g\|_{\bbH_\theta}\|g\|_{\bbH_\theta},
\]
and consequently
\[
 \|g\|_{\bbH_\theta}
 \leq\frac{16}{3}\|\calL g\|_{\bbH_\theta}.
\]
Applied to \(g=\calL^{-1}h\), this gives the first estimate in
\eqref{eq:inverse}.
For the graph norm, the preceding decomposition gives the identity
\[
 \calD_\theta g
 =(\calB_0+2\calK)g-2\calL g.
\]
Together with
\(\|\calB_0\|_{\bbH_\theta\to\bbH_\theta}\leq11/8\) and the
low-norm bound, this identity gives
\[
 \|\calD_\theta g\|_{\bbH_\theta}
 \leq
 \left(
 \frac{22}{3}
 +\frac{32}{3}\|\calK\|_{\bbH_\theta\to\bbH_\theta}
 +2
 \right)
 \|\calL g\|_{\bbH_\theta}.
\]
Combining this estimate with the low-norm bound and the definition of
\(\bbX_\theta\) yields
\[
 \|g\|_{\bbX_\theta}
 \leq
 \left(
 \frac{44}{3}
 +\frac{32}{3}\|\calK\|_{\bbH_\theta\to\bbH_\theta}
 \right)
 \|\calL g\|_{\bbH_\theta}.
\]
Applied to \(\calL^{-1}g\), this gives the second estimate in
\eqref{eq:inverse}.
\end{proof}

\subsection{Coefficient bounds and the affine stationary map}

Having established the linear isomorphism, we next verify that the
nonlinear stationary equation is well defined and analytic on the
affine graph space.  This is the remaining input needed to apply the
implicit function theorem.  Products of perturbations satisfy direct
weighted estimates, whereas terms involving the reference profile
must be grouped so that the zero-slope condition at the origin is
preserved.

For the finite estimates in this subsection and
Lem.~\ref{lem:coefficients}, fix a single constant \(C_\theta\geq1\)
large enough for all of them throughout \(5/4\leq a\leq7/4\).
The proofs give constants depending only on this fixed parameter
interval and the circle; their finite maximum is the one used below.
In particular, \(C_\theta\) is independent of the input functions.

\begin{lemma}[Weighted products]\label{lem:products}
For \(g\in\bbH_\theta\), the bounds
\begin{align}
 &\|g\|_{H^1}+\left\|\frac g{\sin(\theta/2)}\right\|_\infty
 +\|\bmH_{\bbT}g\|_{H^1}
 +|\mathfrak u_1(g)|+|\mathfrak c_a(g)|\notag\\
 &\hspace{12mm}+\|(\calP g)_\theta\|_\infty
 +\left\|\frac{\calP g}{\sin\theta}\right\|_\infty
 +\|\calP g\cot(\theta/2)\|_\infty
 \leq C_\theta\|g\|_{\bbH_\theta}
 \label{eq:basic-bounds}
\end{align}
hold. For \(g_1\in\bbH_\theta\) and \(g_2\in\bbX_\theta\),
the products satisfy
\begin{equation}\label{eq:products}
 \|(\calP g_1)(g_2)_\theta\|_{\bbH_\theta}
 \leq C_\theta\|g_1\|_{\bbH_\theta}\|g_2\|_{\bbX_\theta},
 \qquad
 \|(\bmH_{\bbT}g_1)g_2\|_{\bbH_\theta}
 \leq C_\theta\|g_1\|_{\bbH_\theta}\|g_2\|_{\bbH_\theta}.
\end{equation}
Moreover, \(g\mapsto\calG(a,\widetilde{\bar\Omega}_{3/2}+g)\)
is a real analytic map \(\bbX_\theta\to\bbH_\theta\).
\end{lemma}
\begin{proof}
The basic bounds follow directly from the weighted derivative
control.  For \(0\leq\theta\leq\pi\), the fundamental theorem of
calculus and Cauchy--Schwarz give the estimate
\[
 |g(\theta)|
 \leq
 \left\|\frac{g_\theta}{\sin(\theta/2)}\right\|_{L^2(0,\pi)}
 \left(\int_0^\theta\sin^2(s/2)\idiff s\right)^{1/2}.
\]
The ratio of the second factor to \(\sin(\theta/2)\) is uniformly
bounded on \([0,\pi]\).  Together with oddness, this gives
\[
 \left\|\frac{g}{\sin(\theta/2)}\right\|_\infty
 +\|g\|_{H^1}
 \leq C_\theta\|g\|_{\bbH_\theta}.
\]
The periodic Hilbert transform has Fourier multiplier of modulus one
on the nonzero modes and commutes with differentiation.  Its
\(H^1\) bound, together with the one-dimensional Sobolev embedding,
therefore controls both
\(\|\bmH_{\bbT}g\|_{H^1}\) and
\(\| \bmH_{\bbT}g\|_\infty\) by
\(C_\theta\|g\|_{\bbH_\theta}\).  The definitions of
\(\mathfrak u_1\) and \(\mathfrak c_a\), with \(a\) in the fixed
parameter interval, then give the corresponding bounds for these two
functionals.
It remains to control \(\calP g\) at both endpoints.  The definition
of \(\mathfrak u_1\) gives the vanishing total integral
\[
 \int_0^\pi
 \cos(s/2)
 \left(\frac12(\bmH_{\bbT}g)(s)-\mathfrak u_1(g)\right)\idiff s
 =0.
\]
Consequently, the two representations
\[
 (\calP g)(\theta)
 =
 \frac{1}{\cos(\theta/2)}
 \int_0^\theta
 \cos(s/2)
 \left(\frac12(\bmH_{\bbT}g)(s)-\mathfrak u_1(g)\right)\idiff s
\]
and
\[
 (\calP g)(\theta)
 =
 -\frac{1}{\cos(\theta/2)}
 \int_\theta^\pi
 \cos(s/2)
 \left(\frac12(\bmH_{\bbT}g)(s)-\mathfrak u_1(g)\right)\idiff s
\]
give the endpoint bounds
\[
 |(\calP g)(\theta)|
 \leq
 C_\theta\|g\|_{\bbH_\theta}
 \min\{\theta,\pi-\theta\}.
\]
In particular, the quotient estimate
\[
 \left\|\frac{\calP g}{\sin\theta}\right\|_\infty+\|\calP g\cot(\theta/2)\|_\infty
 \leq C_\theta\|g\|_{\bbH_\theta}
\]
follows.  Finally, the differential identity
\eqref{eq:P-ode}, together with the same quotient estimate and the
\(L^\infty\) bound for \(\bmH_{\bbT}g\), gives
\(
 \|(\calP g)_\theta\|_\infty
 \leq C_\theta\|g\|_{\bbH_\theta}.
\)
These estimates establish \eqref{eq:basic-bounds}.

We now use these pointwise and endpoint bounds to control the
nonlinear products appearing in the stationary equation.
For the first product, the identity
\[
 \partial_\theta\bigl((\calP g_1)(g_2)_\theta\bigr)
 =
 \frac{\calP g_1}{\sin\theta}
 (\calD_\theta g_2)_\theta
 +
 \left(
 (\calP g_1)_\theta
 -\frac{\calP g_1\cos\theta}{\sin\theta}
 \right)(g_2)_\theta
\]
separates the graph derivative of \(g_2\) from the lower-order
weighted derivative.  After division by \(\sin(\theta/2)\), the
first term is controlled by
\(\|\calP g_1/\sin\theta\|_\infty\) and
\(\|\calD_\theta g_2\|_{\bbH_\theta}\), while the second is
controlled by
\[
 \left\|
 (\calP g_1)_\theta
 -\frac{\calP g_1\cos\theta}{\sin\theta}
 \right\|_\infty
\]
and \(\|g_2\|_{\bbH_\theta}\).  The bounds in
\eqref{eq:basic-bounds} therefore give the first estimate in
\eqref{eq:products}.
For the second product, differentiation gives the decomposition
\[
 \partial_\theta\bigl((\bmH_{\bbT}g_1)g_2\bigr)
 =
 (\bmH_{\bbT}g_1)_\theta g_2
 +(\bmH_{\bbT}g_1)(g_2)_\theta .
\]
After division by \(\sin(\theta/2)\), the first term is controlled by
\(\|(\bmH_{\bbT}g_1)_\theta\|_2\) together with
\(\|g_2/\sin(\theta/2)\|_\infty\), whereas the second is controlled by
\(\|\bmH_{\bbT}g_1\|_\infty\) together with
\(\|g_2\|_{\bbH_\theta}\).  The Hilbert-transform and weighted
embedding bounds in \eqref{eq:basic-bounds} then give the second
estimate in \eqref{eq:products}.

It remains to treat the affine terms in which the perturbation is
paired with the reference profile.  These terms require an additional
cancellation at the origin, since neither factor can be estimated
separately in the weighted space.
Evaluation of \eqref{eq:P-ode} at the origin gives the identity
\[
 (\calP g)_\theta(0)
 =\frac12(\bmH_{\bbT}g)(0)-\mathfrak u_1(g),
\]
and hence the representation
\[
 \mathfrak c_a(g)
 =a(\calP g)_\theta(0)-\frac12(\bmH_{\bbT}g)(0).
\]
Together with
\(\widetilde{\bar\Omega}_{3/2}=-\sin\theta\), this yields the
cancellation identity
\begin{align}\label{eq:affine-cancellation}
 &-a(\calP g)(\widetilde{\bar\Omega}_{3/2})_\theta
 +\left(\mathfrak c_a(g)+\frac12\bmH_{\bbT}g\right)
 \widetilde{\bar\Omega}_{3/2}\notag\\
 ={}&-\frac12
 \bigl(\bmH_{\bbT}g-(\bmH_{\bbT}g)(0)\bigr)\sin\theta+
 a\bigl((\calP g)\cos\theta-(\calP g)_\theta(0)\sin\theta\bigr) .
\end{align}
Both terms on the right have zero slope at the origin, as required by
the definition of \(\bbH_\theta\).
For the Hilbert-transform term, differentiation and division by
\(\sin(\theta/2)\) produce terms controlled by
\(\|(\bmH_{\bbT}g)_\theta\|_2\) and by the Hardy estimate
\[
 \left\|
 \frac{\bmH_{\bbT}g-(\bmH_{\bbT}g)(0)}
      {\sin(\theta/2)}
 \right\|_2
 \leq
 C_\theta\|(\bmH_{\bbT}g)_\theta\|_2.
\]
The \(H^1\) boundedness of the Hilbert transform therefore gives the
required \(\bbH_\theta\) bound.
For the \(\calP g\) term, subtraction of the value of
\eqref{eq:P-ode} at the origin gives the identity
\[
 \frac{(\calP g)_\theta-(\calP g)_\theta(0)}
      {\sin(\theta/2)}
 =
 \frac12
 \frac{\bmH_{\bbT}g-(\bmH_{\bbT}g)(0)}
      {\sin(\theta/2)}
 +
 \frac{\calP g}{2\cos(\theta/2)}.
\]
The first term is controlled by the same Hardy estimate, while the
second is controlled by the endpoint estimate for \(\calP g\).  Differentiation of
\((\calP g)\cos\theta-(\calP g)_\theta(0)\sin\theta\) then gives the
corresponding \(\bbH_\theta\) bound.  Consequently,
\eqref{eq:affine-cancellation} defines a bounded map from
\(\bbH_\theta\) to \(\bbH_\theta\).

To finish the analyticity statement, expand the stationary map around
the reference profile.  The decomposition
\[
 \begin{aligned}
 \calG(a,\widetilde{\bar\Omega}_{3/2}+g)
 ={}&
 \calG(a,\widetilde{\bar\Omega}_{3/2})\\
 &-a(\calP\widetilde{\bar\Omega}_{3/2})g_\theta
 +\left(
   \mathfrak c_a(\widetilde{\bar\Omega}_{3/2})
   +\frac12\bmH_{\bbT}\widetilde{\bar\Omega}_{3/2}
  \right)g\\
 &-a(\calP g)(\widetilde{\bar\Omega}_{3/2})_\theta
 +\left(\mathfrak c_a(g)+\frac12\bmH_{\bbT}g\right)
  \widetilde{\bar\Omega}_{3/2}\\
 &-a(\calP g)g_\theta
 +\left(\mathfrak c_a(g)+\frac12\bmH_{\bbT}g\right)g
 \end{aligned}
\]
separates the smooth residual, the linear terms, and the quadratic
terms.  The first linear pair is controlled by smooth multipliers and
\(\calD_\theta g\), the second linear pair by
\eqref{eq:affine-cancellation}, and the quadratic terms by
\eqref{eq:products} together with the bound for
\(\mathfrak c_a(g)\).  Thus
\(g\mapsto\calG(a,\widetilde{\bar\Omega}_{3/2}+g)\) is a quadratic,
and hence real analytic, map from \(\bbX_\theta\) to
\(\bbH_\theta\).
\end{proof}

\subsection{Transversality and the stationary branch}

The linear isomorphism and the analytic stationary map now allow us
to construct the nearby profile family by the implicit function
theorem.  To determine its physical behavior, we must also
track the amplitude scaling rate along the branch.  We therefore
prove a transversality estimate at \(a=3/2\), which determines the
local sign of that rate.

The scaling trace is the functional
\begin{equation}\label{eq:trace-functional}
 \mathfrak l(g)=(\calR g)_\theta(0)=2\mathfrak c_{3/2}(g).
\end{equation}
The basis trace identity \eqref{eq:trace-basis} and the coefficient
bounds \eqref{eq:beta-bounds} in Appendix~\ref{app:rows} give
\begin{equation}\label{eq:trace-bound}
\begin{aligned}
 \|\mathfrak l\|_{\bbH_\theta^*}^2
 &=\sum_{j\geq1}|\mathfrak l(e_j)|^2
 =\sum_{j\geq1}|\beta_j^{\rm tr\text{-}diff}|^2\leq
\underbrace{\sum_{j=1}^4|\beta_j^{\rm tr\text{-}diff}|^2}_{\text{explicit low modes}}
+
\underbrace{\sum_{j\geq5}|\beta_j^{\rm tr\text{-}diff}|^2}_{\text{tail}}\\
 &\leq\left(\frac3{10}\right)^2
 +\left(\frac{31}{210}\right)^2
 +\left(\frac3{140}\right)^2
 +\left(\frac{157}{4620}\right)^2
 +\frac3{1600}
 =\frac{9842363}{85377600}<\left(\frac{17}{50}\right)^2.
 \end{aligned}
\end{equation}
Here the tail follows from
\(\sum_{j=5}^\infty j^{-4}\leq\int_4^\infty s^{-4}\idiff s=1/192\).
The residual of the stationary equation is
\(\calG(a,\widetilde{\bar\Omega}_{3/2})=(a-3/2)e_1/3\).
Differentiating the branch
identity
\(\calG(a,\widetilde{\bar\Omega}_a)=0\) and \(
\bar c_{\omega,a}
=
\mathfrak c_a(\widetilde{\bar\Omega}_a)
\) at \(a=3/2\) gives the tangent and rate identities
\[
 \begin{aligned}
 \calL\left(\left.\partial_a\widetilde{\bar\Omega}_a
                              \right|_{a=3/2}\right)
 &=-\frac13e_1,\\
 \left.\partial_a\bar c_{\omega,a}\right|_{a=3/2}
 &=\left.\partial_a\mathfrak c_a
                 (\widetilde{\bar\Omega}_{3/2})\right|_{a=3/2}
   +\mathfrak c_{3/2}\left(
        \left.\partial_a\widetilde{\bar\Omega}_a\right|_{a=3/2}\right)\\
 &=\frac13+\frac12\mathfrak l\left(
        \left.\partial_a\widetilde{\bar\Omega}_a\right|_{a=3/2}\right).
 \end{aligned}
\]
We therefore define the tangent perturbation and its scaling derivative by
\begin{equation}\label{eq:transversality}
 f_*=-\frac13\calL^{-1}e_1,\qquad
 \kappa_* =\frac13+\frac12\mathfrak l(f_*).
\end{equation}
The inverse estimate in Prop.~\ref{prop:linear} together with the
normalization \(\|e_1\|_{\bbH_\theta}=1\) gives \(\|f_*\|_{\bbH_\theta}\leq16/9\).
Hence \eqref{eq:trace-bound} gives
\[
 \kappa_*
 \geq
 \frac13-\frac12\cdot\frac{136}{225}
 =\frac7{225}>0.
\]
This strict positivity is the required transversality at \(a=3/2\).
It determines the local sign change of the amplitude scaling rate
along the profile branch.

\begin{proposition}[The normalized branch]\label{prop:branch}
For sufficiently small \(\delta_{\rm br},r_{\rm br}>0\), there
is a unique real analytic family
\(\widetilde{\bar\Omega}_a\in\widetilde{\bar\Omega}_{3/2}+\bbX_\theta\)
of zeros of \(\calG(a,\cdot)\) in the radius-\(r_{\rm br}\)
neighborhood of \(\widetilde{\bar\Omega}_{3/2}\). Its perturbation
\(f_a=\widetilde{\bar\Omega}_a-\widetilde{\bar\Omega}_{3/2}\) satisfies
\(\|f_a\|_{\bbX_\theta}\leq C_{\rm br}|a-3/2|\), and its rates obey
\begin{equation}\label{eq:branch-rates}
 \bar c_{\omega,a}=\kappa_*(a-3/2)+O((a-3/2)^2),
 \qquad \bar c_{l,a}=-\frac a6+O(|a-3/2|)<0.
\end{equation}
The constants in these expansions are uniform on the fixed branch
neighborhood. The real representative
\(\bar\Omega_a=\calT^{-1}\widetilde{\bar\Omega}_a\)
is smooth in \((-1,1)\), negative on \((0,1)\), and its zero
extension satisfies \eqref{eq:profile} distributionally.
\end{proposition}
\begin{proof}
Consider the stationary map
\[
 (a,f)\longmapsto
 \calG(a,\widetilde{\bar\Omega}_{3/2}+f)
\]
from a neighborhood of \((3/2,0)\) in
\(\bbR\times\bbX_\theta\) to \(\bbH_\theta\).
Lemma~\ref{lem:products} gives its real analyticity, while
Prop.~\ref{prop:linear} identifies its derivative with respect to
\(f\) at the reference point as the isomorphism \(\calL\).
The Banach implicit function theorem therefore gives a unique real
analytic branch \(f_a\in\bbX_\theta\) with \(f_{3/2}=0\).
The boundedness of the derivative of this analytic branch on a fixed smaller parameter neighborhood gives the estimate
\(
 \|f_a\|_{\bbX_\theta}
 \leq C_{\rm br}|a-3/2|.
\)
The tangent of the branch at the reference parameter is \(f_*\).
The identities in \eqref{eq:transversality} give the first rate
expansion
\[
 \bar c_{\omega,a}
 =\kappa_*(a-3/2)+O((a-3/2)^2).
\]
For the length-scaling rate,  the circle normalization gives the
identity
\(
 \bar c_{l,a}
 =-a\,\mathfrak u_1(\widetilde{\bar\Omega}_a).
\)
The reference value
\(\mathfrak u_1(\widetilde{\bar\Omega}_{3/2})=1/6\), together with
the continuity of \(\mathfrak u_1\) and the branch estimate, gives
the expansion
\[
 \bar c_{l,a}
 =-\frac a6+O(|a-3/2|).
\]
By choosing \(\delta_{\rm br}>0\) sufficiently small, the strict
negativity of \(\bar c_{l,a}\) persists throughout the branch
neighborhood.
The same choice of parameter neighborhood also preserves the
nondegeneracy of the effective transport. 
The bound in \eqref{eq:basic-bounds} and the branch estimate
give the perturbative estimate
\[
 \left\|
 \frac{\calP(\widetilde{\bar\Omega}_a
             -\widetilde{\bar\Omega}_{3/2})}
      {\sin\theta}
 \right\|_\infty
 \leq C_\theta\|f_a\|_{\bbH_\theta}\leq C_\theta C_{\rm br}|a-3/2|.
\]
Shrinking the parameter neighborhood, the identity \((\calP\widetilde{\bar\Omega}_{3/2})/\sin\theta=1/3\) therefore gives the
uniform lower bound
\[
 \frac{\calP\widetilde{\bar\Omega}_a}{\sin\theta}
 \geq\frac16.
\]

It remains to verify the qualitative properties of the real profile:
interior smoothness and sign, regularity through the origin, and the
validity of its zero extension across the support endpoints.
On every compact subinterval of \((0,\pi)\), the lower bound for
\(\calP\widetilde{\bar\Omega}_a/\sin\theta\) keeps the transport
coefficient in the stationary equation uniformly away from zero.
The stationary relation
\[
 a(\calP\widetilde{\bar\Omega}_a)
   (\widetilde{\bar\Omega}_a)_\theta
 =
 \left(
 \mathfrak c_a(\widetilde{\bar\Omega}_a)
 +\frac12\bmH_{\bbT}\widetilde{\bar\Omega}_a
 \right)
 \widetilde{\bar\Omega}_a
\]
can therefore be viewed locally as a nondegenerate first-order
equation for the profile.  Local Sobolev regularity of the periodic
Hilbert transform, together with a cutoff separating the local and
distant parts of its kernel, then bootstraps the profile by one
derivative at each step.  Iteration gives smoothness on
\((0,\pi)\), and hence smoothness of the real profile away from
the origin and the support endpoints. Regularity through the origin
is proved next.
The weighted normalization
\(\widetilde{\bar\Omega}_a(\theta)/\theta\to-1\) gives strict
negativity near the origin.  On every compact interior interval,
the stationary equation is a scalar linear equation for the profile
once its nonlocal coefficient is fixed.  Uniqueness for this equation
prevents a nontrivial profile from acquiring an interior zero.
Consequently, \(\widetilde{\bar\Omega}_a\) is strictly negative on
\((0,\pi)\), and \(\bar\Omega_a\) is strictly negative on
\((0,1)\).

To bootstrap regularity through the origin, we factor out the simple
odd zero and define
the even quotient \(q_a(y)=\bar\Omega_a(y)/y\), with its continuous
value at zero and zero extension outside \([-1,1]\). The graph identity \eqref{eq:weighted-graph} and the relation
\[
 2\cos(\theta/2)(f_a)_{\theta\theta}
 =\frac{(\calD_\theta f_a)_\theta}{\sin(\theta/2)}
   -\cos\theta\,\frac{(f_a)_\theta}{\sin(\theta/2)}
\]
give \(f_a\in H^2\) near zero. Since
\(q_a(y)=\int_0^1(\bar\Omega_a)_y(ry)\idiff r\),
the Hardy averaging bound then places \(q_a\) in \(H^1\) locally. The commutator identity
\begin{equation}\label{eq:origin-commutator}
 \bmH_{\bbR}\bar\Omega_a(y)
 =y\bmH_{\bbR}q_a(y)-\frac1\pi\int_{\bbR}q_a(s)\idiff s
\end{equation}
gives the two identities
\[
 \begin{aligned}
 \bmH_{\bbR}\bar\Omega_a(y)-\bmH_{\bbR}\bar\Omega_a(0)
   &=y\bmH_{\bbR}q_a(y),\\
 \bar U_a(y)&=y\bmH_{\bbR}\bar\Omega_a(0)
              +\int_0^y s\bmH_{\bbR}q_a(s)\idiff s.
 \end{aligned}
\]
Substitute \(\bar\Omega_a=yq_a\) in \eqref{eq:profile} and use
\(\bar c_{\omega,a}=\bar c_{l,a}+(a-1)\bmH_{\bbR}\bar\Omega_a(0)\).
The constant strain terms cancel before division by \(y\), giving
\begin{equation}\label{eq:origin-bootstrap}
 \left(\bar c_{l,a}+a\frac{\bar U_a(y)}y\right)(q_a)_y
 =\left(\bmH_{\bbR}q_a(y)
 -\frac a{y^2}\int_0^y s\bmH_{\bbR}q_a(s)\idiff s\right)q_a(y).
\end{equation}
At the reference profile, the coefficient of \((q_a)_y\) at the
origin has the value
\[
 \bar c_{l,3/2}
 +\frac32(\bar U_{3/2})_y(0)
 =\frac12.
\]
Continuity of the branch and of the corresponding traces therefore
keeps this coefficient uniformly away from zero near \(y=0\) after
the parameter neighborhood is reduced. A cutoff equal to
one near zero splits the Hilbert transform into a local Sobolev
operator and a smooth distant integral. The integral on the right has the regular representation
\[
 \frac1{y^2}\int_0^y s\bmH_{\bbR}q_a(s)\idiff s
 =\int_0^1r\bmH_{\bbR}q_a(ry)\idiff r.
\]
For each nonnegative integer \(\ell\), differentiation brings down
the factor \(r^\ell\); Minkowski's inequality bounds the resulting
local \(L^2\) norm by the integrable factor
\(\int_0^1r^{\ell+1/2}\idiff r=(\ell+3/2)^{-1}\).
Thus this average preserves every local Sobolev class, and
\eqref{eq:origin-bootstrap} first gives \(q_a\in H^2\) locally
and then, by iteration, every \(H^\ell\) class. This also
upgrades the origin slope to the classical value \(-1\).

Finally, we verify that the zero extension satisfies the stationary
equation across the support endpoints.  The profile has zero trace
there, while the endpoint bound
\(\calP\widetilde{\bar\Omega}_a=O(\pi-|\theta|)\) gives the
quadratic vanishing of the real effective velocity,
\[
 \bar c_{l,a}y+a\bar U_a(y)
 =
 \frac a2
 (\calP\widetilde{\bar\Omega}_a)(\theta)\cos(\theta/2)
 =
 O((\pi-|\theta|)^2).
\]
The zero trace and the local \(W^{1,1}\) regularity of the profile
give a locally absolutely continuous zero extension across
\(y=\pm1\).  The preceding vanishing of the effective velocity
eliminates the endpoint flux in the distributional stationary
equation.  Hence the zero extension of \(\bar\Omega_a\) satisfies
\eqref{eq:profile} on the whole real line. Interior smoothness does not assert a smooth
zero extension.
\end{proof}

\section{Nonlinear stability and physical reconstruction}\label{sec:stability}

This section establishes nonlinear stability of the nearby stationary
profile branch and reconstructs the associated physical solutions.
We start with coefficient and multiplier estimates adapted to the
normalized slope and the commuted equation.  These estimates lead to
two complementary energy bounds.  A coercive low-norm estimate
controls the size and decay of the perturbation, while a graph-norm
estimate provides continuation without requiring smallness of the
initial graph norm.  We then develop the corresponding local
well-posedness theory on the maximal graph domain, including
continuous dependence of the solution map.  Together, the local
theory and the two energy estimates yield global stability in
rescaled time.  We conclude by returning to the real line, constructing
smooth compactly supported data within the stability basin, and
integrating the scaling rates to recover the physical lifetime and
the asymptotically self-similar behavior.

We begin by isolating the coefficient structure of the normalized
equation.  The following quantities are chosen so that the endpoint
behavior and the cancellation preserving the normalized origin slope
remain explicit after commuting the equation.
For odd \(H^1(\bbT)\) inputs, define the coefficient families by
\begin{equation}\label{eq:coefficients}
 \begin{gathered}
 B_g=\frac{\calP g}{\sin\theta},\qquad
 C_g=\mathfrak c_a(g)+\frac12\bmH_{\bbT}g,\qquad
 F_g=\frac12\bmH_{\bbT}g-\mathfrak u_1(g),\\
 T_g=C_g-a(B_g+\calD_\theta B_g),\qquad
 \Xi_g=\calD_\theta\bmH_{\bbT}g.
 \end{gathered}
\end{equation}
The quotients have their regular endpoint traces. The equation
\begin{equation}\label{eq:coefficient-ode}
 \calD_\theta B_g+\frac{3\cos\theta-1}{2}B_g=F_g
\end{equation}
follows from \eqref{eq:P-ode}. At the origin it gives
\(B_g(0)=F_g(0)\), \(\calD_\theta B_g(0)=0\),
\(C_g(0)=aB_g(0)\), and \(T_g(0)=0\).
The last two equalities express preservation of the normalized slope.

The even weighted space needed for the nonlocal source is
\begin{equation}\label{eq:even-space}
 \bbY=\left\{g:\ g\text{ is even on }\bbT,\ g(0)=0,\
 \frac{g_\theta}{\sin(\theta/2)}\in L^2\right\},\qquad
 \|g\|_{\bbY}=\left\|\frac{g_\theta}{\sin(\theta/2)}\right\|_2.
\end{equation}
The trace in this definition is the continuous trace supplied by
the weighted derivative.

\begin{lemma}[Coefficient and multiplier estimates]\label{lem:coefficients}
The fixed constant \(C_\theta\) from Lem.~\ref{lem:products} can
be chosen to give the bound
\begin{equation}\label{eq:coefficient-bounds}
 \|B_g\|_{H^1}+\|\calD_\theta B_g\|_{H^1}
 +\|C_g\|_{H^1}+\|T_g\|_{H^1}
 \leq C_\theta\|g\|_{\bbH_\theta}
\end{equation}
for inputs \(g\) in \(\bbH_\theta\). For inputs in \(\bbX_\theta\),
the identity and estimate
\begin{equation}\label{eq:Xi}
 \Xi_g=\bmH_{\bbT}\calD_\theta g
 -(\bmH_{\bbT}\calD_\theta g)(0),\qquad
 \|\Xi_g\|_{\bbY}
 +\left\|\frac{\Xi_g}{\sin(\theta/2)}\right\|_\infty
 \leq C_\theta\|\calD_\theta g\|_{\bbH_\theta}
\end{equation}
hold.
For an even multiplier \(M\in H^1(\bbT)\), the bounds
\begin{equation}\label{eq:multipliers}
 \|Mg\|_{\bbH_\theta}\leq C_\theta\|M\|_{H^1}\|g\|_{\bbH_\theta},
 \qquad
 \|M(g_{\rm ref}+g)\|_{\bbH_\theta}
 \leq C_{\rm ref}\|M\|_{H^1}(1+\|g\|_{\bbH_\theta})
\end{equation}
hold, with the second bound requiring \(M(0)=0\). Here
\(g_{\rm ref}\) is a fixed smooth odd periodic function and
\(C_{\rm ref}\) is a fixed constant depending only on that function.
\end{lemma}
\begin{proof}
The formula for the quotient is
\[
 B_g(\theta)=
 \frac{\int_0^\theta\cos(s/2)F_g(s)\idiff s}
 {2\sin(\theta/2)\cos^2(\theta/2)}.
\]
Near zero this is a smooth multiplier of a first-order Hardy
average, and near \(\pi\) its vanishing total numerator makes it
a smooth multiplier of
\(-(\pi-\theta)^{-2}\int_0^{\pi-\theta}sF_g(\pi-s)\idiff s\).
The averages \(\int_0^1g(rs)\idiff s\) and
\(\int_0^1s g(rs)\idiff s\) preserve \(H^1\). Minkowski's
integral inequality bounds their \(L^2\) operator norms by
\(\int_0^1s^{-1/2}\idiff s=2\) and
\(\int_0^1s^{1/2}\idiff s=2/3\), respectively.
Their differentiated operator norms are at most \(2/3\) and \(2/5\).
A fixed partition of unity between the two
endpoints therefore gives the \(H^1\) bound for \(B_g\).
The coefficient identity \eqref{eq:coefficient-ode} then gives the
\(H^1\) bound for \(\calD_\theta B_g\).  The definitions of \(C_g\)
and \(T_g\), together with the \(H^1\) boundedness of the periodic
Hilbert transform and the already established bounds for \(B_g\) and
\(\calD_\theta B_g\), give the remaining estimates in
\eqref{eq:coefficient-bounds}.

For the differentiated Hilbert term, the commutator is a constant
\[
 (\calD_\theta\bmH_{\bbT}
 -\bmH_{\bbT}\calD_\theta)g
=
 \frac{1}{2\pi}\int_{-\pi}^{\pi}g(s)\sin s\,\idiff s,
\]
as follows by checking the sine modes.  Since
\(\Xi_g=\calD_\theta\bmH_{\bbT}g\) has zero trace at the origin,
evaluation of this commutator there identifies the constant and gives
the representation
\[
 \Xi_g
 =
 \bmH_{\bbT}\calD_\theta g
 -(\bmH_{\bbT}\calD_\theta g)(0).
\]
To estimate its weighted derivative, the basis identity
\[
 \frac{(\bmH_{\bbT}e_j)_\theta}{\sin(\theta/2)}
 =2\cos((j+1/2)\theta)
\]
identifies the relevant weighted derivatives with an orthogonal
half-integer cosine system.  Parseval therefore gives the
\(\bbY\)-norm estimate in \eqref{eq:Xi}.  Since \(\Xi_g(0)=0\),
the fundamental theorem of calculus and Cauchy--Schwarz, with the
same origin weight as in Lem.~\ref{lem:products}, give the quotient
bound
\[
 \left\|\frac{\Xi_g}{\sin(\theta/2)}\right\|_\infty
 \leq C_\theta\|\Xi_g\|_{\bbY}.
\]
Together these estimates give \eqref{eq:Xi}.

We finally establish the multiplier estimates.  Differentiation of the
product gives the identity
\[
 \frac{(Mg)_\theta}{\sin(\theta/2)}
 =
 M\frac{g_\theta}{\sin(\theta/2)}
 +
 M_\theta\frac{g}{\sin(\theta/2)}.
\]
The one-dimensional Sobolev embedding for \(M\), together with the
weighted quotient bound for \(g\), gives the first estimate in
\eqref{eq:multipliers}.

For the affine factor \(g_{\rm ref}+g\), the perturbative part is
controlled in the same way.  The fixed reference function may have a
nonzero slope at the origin, so the product \(Mg_{\rm ref}\) requires
one additional cancellation.  Under the condition \(M(0)=0\), the
Hardy estimate
\[
 \left\|\frac{M}{\sin(\theta/2)}\right\|_2
 \leq C_\theta\|M_\theta\|_2
\]
controls the term containing the reference slope.  The remaining
terms contain only smooth bounded factors from \(g_{\rm ref}\).
This gives the second estimate in \eqref{eq:multipliers}.
\end{proof}

\subsection{The low and graph estimates}
The stability argument relies on two complementary estimates.  A
coercive estimate in the low norm controls the size and decay of the
perturbation.  The graph-norm estimate serves a different purpose.  It
provides continuation on finite rescaled-time intervals and requires
only finiteness, rather than smallness, of the initial graph norm.
We begin with the low-norm estimate and then commute the normalized
equation with \(\calD_\theta\), preserving the zero-slope cancellation
needed to control the graph norm.

Write the circle perturbation and its graph derivative as
\begin{equation}\label{eq:perturbation}
 \begin{gathered}
 f=\widetilde\Omega-\widetilde{\bar\Omega}_a,\qquad
 p=\calD_\theta f,\qquad
 f_0=f(\cdot,0),\quad p_0=\calD_\theta f_0,\\
 f_\tau=\calL_a f+\calN_a(f),
 \end{gathered}
\end{equation}
where \(\calL_a\) is the derivative of \(\calG(a,\cdot)\) at
\(\widetilde{\bar\Omega}_a\), and
\begin{equation}\label{eq:nonlinearity}
 \calN_a(f)=-a(\calP f)f_\theta+C_f f
\end{equation}
defines its quadratic remainder. All ensuing estimates on the
small branch use one fixed constant \(C_{\rm st}\geq1\).
It depends only on \(C_\theta\), the reference profile, and a
fixed branch neighborhood with \(\|f_a\|_{\bbX_\theta}\leq1\).
Choose it as a common bound for the finite family of coefficient,
product, and energy estimates below before choosing the stability
smallness threshold.

We first justify the perturbation of coercivity. The transport part
of \((\calL_a-\calL)f\) is
\(-((a-3/2)\sin\theta/3+a\calP f_a)f_\theta\).
The remaining terms are
\begin{equation}\label{eq:linear-difference}
 \begin{split}
 &-a(\calP f)(f_a)_\theta+C_{f_a}f+C_f f_a+\frac{a-3/2}{3}f\\
 &\qquad +(a-3/2)
 \bigl((\calP f)\cos\theta-(\calP f)_\theta(0)\sin\theta\bigr).
 \end{split}
\end{equation}
They define a bounded map from \(\bbH_\theta\) to itself with
norm at most \(C_{\rm st}(|a-3/2|+\|f_a\|_{\bbX_\theta})\).
For the transport term use the exact quadratic identity
\begin{equation}\label{eq:low-transport}
 \left\langle-Vg_\theta,g\right\rangle_{\bbH_\theta}
 =\frac1{8\pi}\int_{-\pi}^{\pi}
 \bigl(-V\cot(\theta/2)-V_\theta\bigr)
 \frac{|g_\theta|^2}{\sin^2(\theta/2)}\idiff\theta.
\end{equation}
Here \(V\) denotes an odd transport coefficient with zero endpoint
values and bounded \(V_\theta,V\cot(\theta/2)\).
On the core, the integration by parts is
\[
 \begin{aligned}
 4\pi\langle-Vg_\theta,g\rangle_{\bbH_\theta}
 &=-\int_{-\pi}^{\pi}V_\theta
       \frac{|g_\theta|^2}{\sin^2(\theta/2)}\idiff\theta
   -\frac12\int_{-\pi}^{\pi}
       \frac{V}{\sin^2(\theta/2)}
       \partial_\theta(|g_\theta|^2)\idiff\theta\\
 &=\frac12\int_{-\pi}^{\pi}
       \bigl(-V_\theta-V\cot(\theta/2)\bigr)
       \frac{|g_\theta|^2}{\sin^2(\theta/2)}\idiff\theta.
 \end{aligned}
\]
The zero slope and endpoint conditions remove the boundary terms.
The core approximation in Lemma~\ref{lem:domain} extends the identity
to the graph domain. The bounds in
\eqref{eq:basic-bounds} apply to the displayed perturbed velocity.
Shrinking \(\delta_{\rm br}\) therefore yields
\(\langle\calL_a f,f\rangle_{\bbH_\theta}
\leq-\|f\|_{\bbH_\theta}^2/8\).
The same identity with \(V=a\calP f\), and the multiplier estimate
for \(C_f f\), give
\begin{equation}\label{eq:low-energy}
 \frac12\frac{\idiff}{\idiff\tau}\|f\|_{\bbH_\theta}^2
 \leq-\frac18\|f\|_{\bbH_\theta}^2
 +C_{\rm st}\|f\|_{\bbH_\theta}^3.
\end{equation}

The low estimate does not by itself control the graph norm.  We
therefore commute the normalized equation with \(\calD_\theta\),
using the coefficient structure above to preserve the zero-slope
grouping. Note that \(\calD_\theta C_{\widetilde\Omega}
=\Xi_{\widetilde\Omega}/2\) and the product rule give
\[
 \begin{aligned}
 \calD_\theta\bigl(-aB_{\widetilde\Omega}
                    \calD_\theta\widetilde\Omega
                     +C_{\widetilde\Omega}\widetilde\Omega\bigr)
 ={}&-aB_{\widetilde\Omega}\calD_\theta^2\widetilde\Omega\\
 &+(C_{\widetilde\Omega}-a\calD_\theta B_{\widetilde\Omega})
                      \calD_\theta\widetilde\Omega
   +\frac12\widetilde\Omega\Xi_{\widetilde\Omega}.
 \end{aligned}
\]
Moving one multiple of \(aB_{\widetilde\Omega}
\calD_\theta\widetilde\Omega\) into the principal transport
produces the zero-trace coefficient \(T_{\widetilde\Omega}\)
and the equation
\begin{equation}\label{eq:full-commutation}
 \partial_\tau(\calD_\theta\widetilde\Omega)
 +aB_{\widetilde\Omega}(\calD_\theta-\Id)
       \calD_\theta\widetilde\Omega
 =T_{\widetilde\Omega}\calD_\theta\widetilde\Omega
 +\frac12\widetilde\Omega\Xi_{\widetilde\Omega}.
\end{equation}
Define the fixed profile quantity
\(\Delta_a=(\calD_\theta-\Id)\calD_\theta\widetilde{\bar\Omega}_a\).
The stationary version of \eqref{eq:full-commutation} gives the identity
\begin{equation}\label{eq:fixed-graph}
 aB_{\widetilde{\bar\Omega}_a}\Delta_a
 =T_{\widetilde{\bar\Omega}_a}\calD_\theta\widetilde{\bar\Omega}_a
 +\frac12\widetilde{\bar\Omega}_a\Xi_{\widetilde{\bar\Omega}_a}.
\end{equation}
Both even factors on the right vanish at zero, while the odd
factors are fixed smooth functions plus bounded \(\bbH_\theta\)
perturbations. Lem.~\ref{lem:coefficients} and
\(B_{\widetilde{\bar\Omega}_a}\geq1/6\) therefore imply
\(\|\Delta_a\|_{\bbH_\theta}\leq C_{\rm st}\).

One further fixed graph derivative will be needed for continuous
dependence. The exact formulas
\begin{equation}\label{eq:DT}
 \begin{split}
 T_g&=\frac{1-a}{2}
       \bigl(\bmH_{\bbT}g-(\bmH_{\bbT}g)(0)\bigr)
       -3a\sin^2(\theta/2)B_g,\\
 \calD_\theta T_g&=\frac{1-a}{2}\Xi_g
       -\frac{3a}{2}\sin^2\theta\,B_g
       -3a\sin^2(\theta/2)\calD_\theta B_g
 \end{split}
\end{equation}
follow from \eqref{eq:coefficient-ode}. They put
\(\calD_\theta T_g\) in \(H^1\), with zero trace, for inputs
\(g\) in \(\bbX_\theta\), and also for a fixed smooth affine part.
The identity defining \(\Delta_a\) first gives
\(\calD_\theta^2 f_a\in\bbH_\theta\). Applying \eqref{eq:Xi}
to this derivative controls
\(\calD_\theta\Xi_{\widetilde{\bar\Omega}_a}\) in \(\bbY\).
Differentiate \eqref{eq:fixed-graph} to obtain
\begin{equation}\label{eq:fixed-next-graph}
 \begin{split}
 aB_{\widetilde{\bar\Omega}_a}\calD_\theta\Delta_a
 ={}&(\calD_\theta T_{\widetilde{\bar\Omega}_a})
                  \calD_\theta\widetilde{\bar\Omega}_a
 +T_{\widetilde{\bar\Omega}_a}\calD_\theta^2\widetilde{\bar\Omega}_a\\
 &+\frac12(\calD_\theta\widetilde{\bar\Omega}_a)\Xi_{\widetilde{\bar\Omega}_a}
 +\frac12\widetilde{\bar\Omega}_a\calD_\theta\Xi_{\widetilde{\bar\Omega}_a}
 -a(\calD_\theta B_{\widetilde{\bar\Omega}_a})\Delta_a .
 \end{split}
\end{equation}
Every fixed nonzero-slope odd factor is multiplied by an even
factor with zero trace. The multiplier lemma and division by the
positive \(H^1\) function \(B_{\widetilde{\bar\Omega}_a}\) thus give
\(\|\calD_\theta\Delta_a\|_{\bbH_\theta}\leq C_{\rm st}\).
No smoothness of the zero extension at the endpoints is used.

Subtracting the stationary equation from \eqref{eq:full-commutation}
gives the grouped perturbation equation
\begin{equation}\label{eq:commuted}
 p_\tau+aB_{\widetilde\Omega}(\calD_\theta-\Id)p
 =T_{\widetilde\Omega}p+\frac12\widetilde\Omega\Xi_f+E_f,
 \qquad
 E_f=T_f\calD_\theta\widetilde{\bar\Omega}_a
       -aB_f\Delta_a+\frac12 f\Xi_{\widetilde{\bar\Omega}_a}.
\end{equation}
Its source bounds are
\begin{equation}\label{eq:high-sources}
 \|E_f\|_{\bbH_\theta}\leq C_{\rm st}\|f\|_{\bbH_\theta},
 \qquad
 \|\widetilde\Omega\Xi_f\|_{\bbH_\theta}
 \leq C_{\rm st}(1+\|f\|_{\bbH_\theta})\|p\|_{\bbH_\theta}.
\end{equation}
For an even multiplier \(M\) with the bounds in
\eqref{eq:coefficient-bounds}, integration by parts gives
\begin{equation}\label{eq:high-transport}
 \begin{split}
 \langle M(\calD_\theta-\Id)p,p\rangle_{\bbH_\theta}
 =\frac1{4\pi}\int_{-\pi}^{\pi}
 \frac{\bigl(\frac12\calD_\theta M+
       (\cos\theta-\frac12)M\bigr)|p_\theta|^2
       -M_\theta p p_\theta}{\sin^2(\theta/2)}\idiff\theta.
 \end{split}
\end{equation}
The last term is bounded by
\(\|M_\theta\|_2\|p/\sin(\theta/2)\|_\infty
\|p_\theta/\sin(\theta/2)\|_2\).
Thus its right side extends continuously to \(p\in\bbH_\theta\),
even when \(\calD_\theta p\) is not in that space.
Equations \eqref{eq:commuted}--\eqref{eq:high-transport} imply
\begin{equation}\label{eq:high-energy}
 \frac12\frac{\idiff}{\idiff\tau}\|p\|_{\bbH_\theta}^2
 \leq C_{\rm st}(1+\|f\|_{\bbH_\theta})
       \bigl(\|p\|_{\bbH_\theta}^2+\|f\|_{\bbH_\theta}^2\bigr).
\end{equation}

\subsection{Strong solutions on the graph domain}
The preceding estimates are a priori bounds and are initially
justified only for compatible smooth approximations.  To use them in
the nonlinear stability argument, we need a local well-posedness
theory on the maximal graph domain, with an existence time controlled
by the graph norm but without any smallness assumption on that norm.
We therefore construct strong graph-domain solutions, establish their
continuation criterion, and prove continuous dependence in the graph
topology.

\begin{lemma}[Local strong circle solutions]\label{lem:circle-local}
For every \(f_0\in\bbX_\theta\), the perturbation equation has a
unique maximal solution
\[
 f\in C([0,\tau_*);\bbX_\theta)
       \cap C^1([0,\tau_*);\bbH_\theta).
\]
On each graph ball of radius \(R>0\), a common existence time
\(\tau_{\rm loc}(R)>0\) can be chosen, and the solution map is
continuous in this displayed topology. A finite \(\tau_*\) is
impossible if the graph norm stays bounded up to \(\tau_*\).
\end{lemma}
\begin{proof}
We first solve equations with a frozen velocity by characteristics,
then iterate them with uniform graph bounds. Convergence in the low
norm and weak convergence of graph derivatives identify the limit.
A weighted characteristic formula recovers strong graph time traces.
Finally, smoothing by one additional graph order compensates for
the derivative lost in differences and proves graph continuous
dependence. These steps use a local time determined by the first
graph norm alone.

Define the bilinear evolution map by
\begin{equation}\label{eq:frozen}
 \calB(g_1,g_2)=-a(\calP g_1)(g_2)_\theta+C_{g_1}g_2.
\end{equation}
For \(g_1\) continuous in time in a bounded subset of
\(-\sin\theta+\bbX_\theta\), the velocity
\(a\calP g_1=aB_{g_1}\sin\theta\) is Lipschitz and fixes
\(0,\pm\pi\). The relation \(C_{g_1}(0)=aB_{g_1}(0)\)
preserves the prescribed slope of a solution of
\((g_2)_\tau=\calB(g_1,g_2)\).
Approximate \(g_1+\sin\theta\) and \(g_2(0)+\sin\theta\)
in their graph norms by time-smooth finite-core functions. Restoring
the fixed \(-\sin\theta\) part gives smooth coefficients and data
with exactly the same affine constraint. Characteristics give
smooth linear solutions. Write \(g_{1,N}\) for the restored
coefficient approximations and \(g_{2,N}\) for their evolved
solutions with approximated initial data; \(N,N'\) denote
positive integer approximation indices. Their commuted equation is
\[
 \partial_\tau\calD_\theta g_2
 +aB_{g_1}(\calD_\theta-\Id)\calD_\theta g_2
 =T_{g_1}\calD_\theta g_2+\frac12 g_2\Xi_{g_1}.
\]
Subtract the fixed smooth reference in its low and high equations.
The source for \(\calD_\theta(g_2+\sin\theta)\) is the sum
\[
 T_{g_1}\calD_\theta(g_2+\sin\theta)
 -T_{g_1}\calD_\theta\sin\theta
 +aB_{g_1}(\calD_\theta-\Id)\calD_\theta\sin\theta
 +\frac12 g_2\Xi_{g_1}.
\]
Each grouped term belongs to \(\bbH_\theta\):
\(T_{g_1}\) vanishes at zero,
\((\calD_\theta-\Id)\calD_\theta\sin\theta\) has zero slope,
and \(\Xi_{g_1}\) has the weighted zero trace in \eqref{eq:Xi}.
The multiplier and transport estimates give uniform graph bounds
on a time interval depending only on the common graph bound.
The difference of two approximations has forcing
\(\calB(g_{1,N}-g_{1,N'},g_{2,N'})\), and the identity
\begin{equation}\label{eq:difference-grouping}
 \calB(g_1,g_2)=(C_{g_1}-aB_{g_1})g_2
                +aB_{g_1}(g_2-\calD_\theta g_2)
\end{equation}
places this forcing in \(\bbH_\theta\): the first even factor
vanishes at zero, and the second odd factor has zero slope.
The low energy estimate makes the approximations Cauchy in
\(C_\tau\bbH_\theta\). Their graph derivatives remain bounded.
The limit identification and strong-continuity argument below
therefore give a unique affine graph solution of the frozen equation.

Having constructed the frozen evolution, we now close the nonlinear iteration.
For the nonlinear problem, start with the constant-in-time
\(\widetilde\Omega^{(0)}=\widetilde{\bar\Omega}_a+f_0\) and solve
\begin{equation}\label{eq:iteration}
 \partial_\tau\widetilde\Omega^{(n+1)}
 =\calB(\widetilde\Omega^{(n)},\widetilde\Omega^{(n+1)}),
 \qquad \widetilde\Omega^{(n+1)}(0)=\widetilde{\bar\Omega}_a+f_0.
\end{equation}
Here the sequence index \(n\) labels the nonlinear iteration.
Let \(f^{(n)}=\widetilde\Omega^{(n)}-\widetilde{\bar\Omega}_a\)
and \(p^{(n)}=\calD_\theta f^{(n)}\).
For \(\|f^{(n)}\|_{\bbX_\theta}\leq R\), the low equation
and the exact high equation
\begin{equation}\label{eq:iteration-high}
 \begin{split}
 \partial_\tau p^{(n+1)}
 +aB_{\widetilde\Omega^{(n)}}(\calD_\theta-\Id)p^{(n+1)}
 ={}&T_{\widetilde\Omega^{(n)}}p^{(n+1)}
 +T_{f^{(n)}}\calD_\theta\widetilde{\bar\Omega}_a-aB_{f^{(n)}}\Delta_a\\
 &+\frac12\widetilde{\bar\Omega}_a\Xi_{f^{(n)}}
 +\frac12 f^{(n+1)}\Xi_{\widetilde\Omega^{(n)}}
 \end{split}
\end{equation}
give
\[
 \frac{\idiff}{\idiff\tau}
 \bigl(\|f^{(n+1)}\|_{\bbH_\theta}^2+
       \|p^{(n+1)}\|_{\bbH_\theta}^2\bigr)
 \leq C_{\rm loc}(R)
       \bigl(1+\|f^{(n+1)}\|_{\bbX_\theta}^2+
                    \|f^{(n)}\|_{\bbX_\theta}^2\bigr).
\]
Here \(C_{\rm loc}(R)\) is a fixed constant for this and the
difference estimates on that ball, depending only on \(R,C_{\rm st}\).
Taking \(R>2\|f_0\|_{\bbX_\theta}+1\) and a sufficiently
short interval closes the uniform bound by induction and Gronwall.
Differences of consecutive iterates satisfy the same low transport
estimate with the previous-iterate forcing in
\eqref{eq:difference-grouping}. More precisely, for \(n\geq1\),
the common initial data and Gronwall's inequality give
\[
 \|f^{(n+1)}(\tau)-f^{(n)}(\tau)\|_{\bbH_\theta}
 \leq C_{\rm loc}(R)\int_0^\tau
       \|f^{(n)}(s)-f^{(n-1)}(s)\|_{\bbH_\theta}\idiff s.
\]
After increasing the fixed constant on the local interval, iteration
of this inequality yields
\[
 \|f^{(n+1)}-f^{(n)}\|_{C([0,\tau];\bbH_\theta)}
 \leq\frac{(C_{\rm loc}(R)\tau)^n}{n!}
       \|f^{(1)}-f^{(0)}\|_{C([0,\tau];\bbH_\theta)}.
\]
Thus \(f^{(n)}\) converges in
\(C_\tau\bbH_\theta\), and \(p^{(n)}\) converges along a
subsequence weakly star in \(L^\infty_\tau\bbH_\theta\).

To identify the graph limit, test
\(\langle p^{(n)},g\rangle_{\bbH_\theta}
=\langle f^{(n)},\calD_\theta^*g\rangle_{\bbH_\theta}\)
against time test functions, with
\(g\in\operatorname{Dom}(\calD_\theta^*)\). The limit satisfies
the same identity. The double-adjoint characterization of a closed
densely defined operator puts \(f\) in \(\bbX_\theta\) and gives
the identity \(p=\calD_\theta f\) almost everywhere.
Boundedness of the linear maps
\(g\mapsto(B_g,\calD_\theta B_g,C_g,T_g)\) sends the low-norm
convergence to strong \(H^1\) convergence of the coefficient fields.
The bounded linear map
\(p\mapsto\bmH_{\bbT}p-(\bmH_{\bbT}p)(0)\) is weak-to-weak
continuous into \(\bbY\). The multiplier estimates allow its
product with strongly convergent affine factors. Testing the
transport term on the finite core moves its derivative onto the
test function and uses its zero endpoint velocity. Thus both
grouped equations pass to the limit, with their initial traces.

Here is the strong-time-continuity step. For a high equation of the form
\(p_\tau+aB_{g_1}(\calD_\theta-\Id)p=S\), division of its derivative
by \(\sin(\theta/2)\) gives
\begin{equation}\label{eq:weighted-transport}
 \begin{split}
 \partial_\tau\left(\frac{p_\theta}{\sin(\theta/2)}\right)
 +aB_{g_1}\sin\theta\,\partial_\theta
       \left(\frac{p_\theta}{\sin(\theta/2)}\right)
 ={}&\frac{S_\theta}{\sin(\theta/2)}
       +a(B_{g_1})_\theta\frac p{\sin(\theta/2)}\\
 &-a\left(\calD_\theta B_{g_1}
                 +\frac{3\cos\theta-1}{2}B_{g_1}\right)
       \frac{p_\theta}{\sin(\theta/2)} .
 \end{split}
\end{equation}
Here \(S\) denotes the grouped odd source, either in the frozen
equation or in \eqref{eq:commuted}. It is strongly measurable in
\(L^\infty_\tau\bbH_\theta\); all factors belong to separable
spaces and are connected by bounded linear or bilinear maps.
The right side is therefore in \(L^\infty_\tau L^2_\theta\).
The velocity lies in \(C_\tau W^{1,\infty}_\theta\).
Its characteristic flow fixes \(0,\pm\pi\), and its Jacobian is
the exponential of the integral of the velocity derivative.
Composition by this flow is strongly continuous on \(L^2\),
first for smooth functions and then by density, with norm at most
the square root of the maximal inverse Jacobian. The characteristic
formula for \eqref{eq:weighted-transport} is consequently a continuous
\(L^2\) initial term plus a Bochner integral of an \(L^1_\tau L^2\)
source.
For completeness, choose a fixed smooth even nonnegative kernel
\(\psi\), supported in \((-1,1)\), with integral one, and denote
convolution in \(\theta\) with \(N\psi(N\theta)\) by
\(\calJ_N^\theta\).
Extend a Lipschitz velocity \(V\) that vanishes at the endpoints
by zero. The commutator identity
\begin{equation}\label{eq:transport-commutator}
 [V\partial_\theta,\calJ_N^\theta]g
 =\int N^2\psi'(N(\theta-s))(V(\theta)-V(s))g(s)\idiff s
   +\calJ_N^\theta(V' g)
\end{equation}
gives an \(L^2\) operator norm at most
\((1+\int|s\psi'(s)|\idiff s)\|V'\|_\infty\).
For smooth \(g\) the commutator tends to zero; density and this
uniform bound give convergence for every \(L^2\) input, also in
the time-integrated norm. Mollification and the energy identity
therefore give uniqueness of the weak transport solution. Its
characteristic representation is the one above. The vanishing
endpoint velocity removes the boundary flux.
We obtain \(p\in C_\tau\bbH_\theta\). Closedness now extends
\(\calD_\theta f=p\) to every time, and the low equation gives
\(f_\tau\in C_\tau\bbH_\theta\).

It remains to establish continuous dependence in the graph topology.
The same grouped difference estimate proves uniqueness and
\begin{equation}\label{eq:low-dependence}
 \|f-g\|_{C([0,\tau_{\rm loc}(R)];\bbH_\theta)}
 \leq C_{\rm loc}(R)\|f_0-g_0\|_{\bbH_\theta}
\end{equation}
for two solutions in the graph ball. To upgrade this to graph
continuity, construct a compatible smoother in the coordinate
\(\zeta\) of \eqref{eq:graph-coordinate}. Let \(\calJ_N^\zeta\)
denote convolution in
\(\zeta\) with \(N\psi(N\zeta)\). Choose a fixed smooth
compatibility corrector \(\psi_{\rm c}\), compactly supported
in \(0<\theta<\pi\), with
\(\int_0^\pi\sin(s/2)\psi_{\rm c}(s)\idiff s=1\).
Define \(\calS_N g\) by reconstructing the weighted derivative
\begin{equation}\label{eq:smoother}
 \frac{(\calS_N g)_\theta}{\sin(\theta/2)}
 =\calJ_N^\zeta\left(\frac{g_\theta}{\sin(\theta/2)}\right)
 -\psi_{\rm c}\int_0^\pi\sin\frac s2\,
     \calJ_N^\zeta\left(\frac{g_\theta}{\sin(\theta/2)}\right)(s)\idiff s.
\end{equation}
The convolution is evaluated in the logarithmic coordinate before
pullback to \(s\). The inequality
\(\operatorname{sech}(\zeta-r)\leq e^{|r|}
\operatorname{sech}\zeta\), the graph identification, and symmetry
of \(\psi\) give
\begin{equation}\label{eq:smoothing}
 \|\calS_N g-g\|_{\bbX_\theta}\to0,\qquad
 N\|\calS_N g-g\|_{\bbH_\theta}\to0,\qquad
 \|\calS_N g\|_{\bbX_{\theta,2}}\leq C_{\rm sm}N\|g\|_{\bbX_\theta}.
\end{equation}
The fixed constant \(C_{\rm sm}\) depends only on the two chosen
kernels. The second limit holds first on the weighted \(H^2\)
core by cancellation of the first mollifier moment, and then on
weighted \(H^1\) by density and the uniform first-order error bound.
The compatibility correction has the same estimates.

For data one graph order smoother, apply \(\calD_\theta\) to
the high equation. Its exact form is
\begin{equation}\label{eq:second-graph}
 \partial_\tau(\calD_\theta p)
 +aB_{\widetilde\Omega}(\calD_\theta-\Id)\calD_\theta p
 =\calD_\theta S
 -a(\calD_\theta B_{\widetilde\Omega})(\calD_\theta-\Id)p.
\end{equation}
The highest nonlocal term is
\(\calD_\theta\Xi_f
=\bmH_{\bbT}\calD_\theta p-(\bmH_{\bbT}\calD_\theta p)(0)\).
For the other coefficients use \eqref{eq:coefficient-ode} and
\eqref{eq:DT}; differentiation of the former expresses
\(\calD_\theta^2B_g\) through \(\Xi_g,B_g,\calD_\theta B_g\)
and a fixed smooth multiplier. Every affine odd factor still
has a zero-trace even multiplier. Together with
\eqref{eq:fixed-next-graph}, this bounds the right side of
\eqref{eq:second-graph} by
\(C_{\rm loc}(R)(1+\|\calD_\theta p\|_{\bbH_\theta})\).
To justify this calculation without assuming additional persistence,
run the same frozen iteration with initial data in \(\bbX_{\theta,2}\).
On the already bounded \(\bbX_\theta\) ball, its differentiated
equation gives
\begin{equation}\label{eq:second-iteration}
 \frac{\idiff}{\idiff\tau}
       \|\calD_\theta^2 f^{(n+1)}\|_{\bbH_\theta}^2
 \leq C_{\rm loc}(R)\left(1+
       \|\calD_\theta^2 f^{(n+1)}\|_{\bbH_\theta}^2+
       \|\calD_\theta^2 f^{(n)}\|_{\bbH_\theta}^2\right).
\end{equation}
The derivative of \(\Xi_{f^{(n)}}\) is the only term involving
the highest derivative of the previous iterate; its coefficient
is bounded on that graph ball. Induction and Gronwall on a time
interval depending only on \(R\) bound these norms by a fixed
multiple of \(1+\|\calD_\theta^2 f_0\|_{\bbH_\theta}\), regardless
of its size. Finite-core approximation, followed by weak graph
identification as above, passes this bound to the solution.
The weighted transport argument applied at this extra level gives
the strong initial trace and strong time continuity. Restarting
on finitely many such intervals covers the common local interval.
The transport identity for the resulting solution now yields
\begin{equation}\label{eq:second-bound}
 \|f\|_{C([0,\tau_{\rm loc}(R)];\bbX_{\theta,2})}
 \leq C_{\rm loc}(R)(1+\|f_0\|_{\bbX_{\theta,2}}).
\end{equation}

Subtracting the commuted equations for two solutions isolates the
only lost derivative in
\(aB_{f-g}(\calD_\theta-\Id)\calD_\theta g\).
The same multiplier estimates, \eqref{eq:low-dependence}, and
Gronwall give
\begin{equation}\label{eq:graph-dependence}
 \begin{split}
 \|f-g\|_{C([0,\tau_{\rm loc}(R)];\bbX_\theta)}
 \leq C_{\rm loc}(R)\bigl(&\|f_0-g_0\|_{\bbX_\theta}\\
 &+(1+\|g\|_{L^1(0,\tau_{\rm loc}(R);\bbX_{\theta,2})})
       \|f_0-g_0\|_{\bbH_\theta}\bigr).
 \end{split}
\end{equation}
Choose \(g_0=\calS_N f_0\). The possible factor \(N\) from
\eqref{eq:second-bound} is absorbed by the second limit in
\eqref{eq:smoothing}. This proves convergence in the graph norm
of the smoothed solutions. For initial data converging to \(f_0\),
compare each corresponding solution directly with the one solution
issuing from \(\calS_N f_0\). First fix \(N\) and take the
upper limit along the data sequence in \eqref{eq:graph-dependence};
then let \(N\) tend to infinity. The resulting upper bound is
a fixed multiple of
\[
 \|f_0-\calS_N f_0\|_{\bbX_\theta}
 +(1+N)\|f_0-\calS_N f_0\|_{\bbH_\theta},
\]
which tends to zero by \eqref{eq:smoothing}. This proves graph continuous
dependence; continuity of the low equation gives the \(C^1\)
part. Finally the construction time depends only on the graph
bound, so restarting just before a finite endpoint proves
continuation.
\end{proof}

We now combine the local theory with the low- and graph-norm estimates to obtain global stability in rescaled time.
Take \(\varepsilon_{\rm st}\leq1/(16C_{\rm st})\). While
\(\|f(\tau)\|_{\bbH_\theta}<\varepsilon_{\rm st}\),
\eqref{eq:low-energy} gives
\[
 \frac12\frac{\idiff}{\idiff\tau}\|f\|_{\bbH_\theta}^2
 \leq-\frac1{16}\|f\|_{\bbH_\theta}^2.
\]
Continuity closes this bootstrap for all times of existence.
Transferring by the isometry \(\calT\), we obtain the decay estimate
\begin{equation}\label{eq:stability}
 \|\Omega(\tau)-\bar\Omega_a\|_{\bbH_y}
 \leq e^{-\gamma_{\rm st}\tau}
             \|\omega_0-\bar\Omega_a\|_{\bbH_y},
 \qquad \gamma_{\rm st}=\frac1{16},
\end{equation}
for every datum satisfying \eqref{eq:basin}.
If the squared graph energy is denoted by
\begin{equation}\label{eq:graph-energy}
 \mathfrak E_{\rm gr}(\tau)=
 \|f(\tau)\|_{\bbH_\theta}^2+\|p(\tau)\|_{\bbH_\theta}^2,
\end{equation}
then \eqref{eq:low-energy} and \eqref{eq:high-energy} give a
fixed constant \(C_{\rm gr}\), depending only on \(C_{\rm st}\),
with the growth bound
\begin{equation}\label{eq:graph-global}
 \mathfrak E_{\rm gr}(\tau)\leq\mathfrak E_{\rm gr}(0)
 \exp\left(C_{\rm gr}\tau+
       \frac{C_{\rm gr}}{\gamma_{\rm st}}\|f_0\|_{\bbH_\theta}\right).
\end{equation}
Thus the graph norm is finite on every finite rescaled interval.
Only the low norm in \eqref{eq:basin} is required to be small;
Lemma~\ref{lem:circle-local} then gives global existence with
\[
 \Omega-\bar\Omega_a\in C([0,\infty);\bbX_y)
                     \cap C^1([0,\infty);\bbH_y).
\]

The inverse circle map also gives the distributional real rescaled
solution. Indeed, the identity
\[
 \int_{-1}^1|\Omega_y(y,\tau)|\idiff y
 =\frac12\int_{-\pi}^{\pi}
             |\widetilde\Omega_\theta(\theta,\tau)|\idiff\theta
 <\infty
\]
and the zero endpoint traces justify zero extension in \(W^{1,1}\).
The effective transport velocity vanishes at the support endpoints,
so integration by parts introduces no boundary source. The Hilbert
conjugacy transfers the circle equation to \eqref{eq:rescaled},
and the isometry transfers its strong time topology.

\subsection{Smooth data and physical time}

The stationary profile is smooth in the support interior, but its zero
extension need not be smooth at the support endpoints.  To obtain the
smooth compactly supported physical data required in
Thm.~\ref{thm:profiles}, we first cut off the profile near its
endpoints and show that the resulting data remain inside the low-norm
stability basin.  We then translate rescaled convergence into
asymptotics for the scaling factors and physical time.  Finally, a
real-line continuation argument identifies the reconstruction with
the classical solution issued from the smooth data and establishes
finite-time blowup or global smoothness according to the sign of the
amplitude rate.

Choose a fixed nonnegative smooth cutoff on the interval and obtain
an even family \(\chi_\sigma\), for \(0<\sigma<\pi/4\), satisfying
\begin{equation}\label{eq:cutoff}
 \chi_\sigma=1\ (|\theta|\leq\pi-2\sigma),\qquad
 \chi_\sigma=0\ (|\theta|\geq\pi-\sigma),\qquad
 0<\chi_\sigma\leq1\ (|\theta|<\pi-\sigma),\qquad
 |\chi_\sigma'|\leq C_\chi/\sigma.
\end{equation}
Here \(C_\chi\) is the derivative bound of the fixed cutoff.
Set \(\widetilde\Omega_{0,\sigma}=\chi_\sigma\widetilde{\bar\Omega}_a\)
and \(\omega_{0,\sigma}=\calT^{-1}\widetilde\Omega_{0,\sigma}\),
with zero extension. These functions are smooth, odd, and have
normalized slope \(-1\); the physical data are supported strictly
inside \((-1,1)\). The branch sign proves strict negativity
inside their positive support interval.
The derivative of the cutoff error splits as
\[
 \partial_\theta(\widetilde\Omega_{0,\sigma}
                    -\widetilde{\bar\Omega}_a)
 =(\chi_\sigma-1)\partial_\theta\widetilde{\bar\Omega}_a
     +\chi_\sigma'\widetilde{\bar\Omega}_a.
\]
The weight \(1/\sin^2(\theta/2)\) is bounded near the two
endpoints. The first term therefore tends to zero in the weighted
derivative norm by absolute continuity of the integral. For the
second term, the endpoint trace and Cauchy--Schwarz give
\[
 \sigma^{-2}\int_0^{2\sigma}
 |\widetilde{\bar\Omega}_a(\pi-r)|^2\idiff r
 \leq2\int_0^{2\sigma}
 |(\widetilde{\bar\Omega}_a)_\theta(\pi-r)|^2\idiff r\longrightarrow0.
\]
Oddness gives the same estimate at the negative endpoint.
These estimates give the limit
\(\|\widetilde\Omega_{0,\sigma}-\widetilde{\bar\Omega}_a\|_{\bbH_\theta}
\to0\). The difference belongs to \(\bbX_\theta\): subtract the
fixed \(-\sin\theta\) part, and multiply its graph perturbation
by the fixed smooth multiplier \(\chi_\sigma-1\), which is constant
near zero. Its graph norm need only be finite for each \(\sigma\).
Thus these data enter the stability basin.

We next translate the rescaled convergence into asymptotics for the physical scaling factors and physical time.
For a general solution in that basin, boundedness of the linear
rate functionals and \eqref{eq:stability} imply
\begin{equation}\label{eq:rate-errors}
 e_\omega(\tau):=c_\omega(\tau)-\bar c_{\omega,a},\qquad
 e_l(\tau):=c_l(\tau)-\bar c_{l,a},\qquad
 |e_\omega(\tau)|+|e_l(\tau)|
 \leq C_{\rm rec}e^{-\gamma_{\rm st}\tau}\|f_0\|_{\bbH_\theta}.
\end{equation}
The fixed constant \(C_{\rm rec}\) is a common operator norm bound
for the two rate functionals on the branch interval. Their integrals
\(I_\omega=\int_0^\infty e_\omega(s)\idiff s\) and
\(I_l=\int_0^\infty e_l(s)\idiff s\) are finite and yield the asymptotics
\begin{equation}\label{eq:scale-asymptotics}
 C_\omega=e^{-I_\omega}e^{-\bar c_{\omega,a}\tau}
 (1+O(e^{-\gamma_{\rm st}\tau})),\qquad
 C_l=e^{I_l}e^{\bar c_{l,a}\tau}
 (1+O(e^{-\gamma_{\rm st}\tau})).
\end{equation}
The constants in these asymptotic bounds may depend on the chosen
parameter and initial perturbation, but are uniform in \(\tau\).
For a negative amplitude rate, define \(T_{\rm phys}=t(\infty)\).
Integration of the first reciprocal scaling factor gives
\begin{equation}\label{eq:physical-lifetime}
 T_{\rm phys}-t(\tau)=
 \frac{e^{I_\omega}}{|\bar c_{\omega,a}|}e^{\bar c_{\omega,a}\tau}
 (1+O(e^{-\gamma_{\rm st}\tau})),\qquad
 (T_{\rm phys}-t(\tau))C_\omega(\tau)\longrightarrow
 \frac1{|\bar c_{\omega,a}|}.
\end{equation}
Uniform convergence of \(\Omega\) to the nonzero profile follows
from the weighted embedding. The reconstruction
\(\omega(x,t(\tau))=C_\omega(\tau)\Omega(C_l(\tau)x,\tau)\)
therefore has unbounded amplitude at the finite time \(T_{\rm phys}\).
For a positive amplitude rate, integration of
\eqref{eq:scale-asymptotics} gives
\[
 t(\tau)=\frac{e^{I_\omega}}{\bar c_{\omega,a}}
              e^{\bar c_{\omega,a}\tau}(1+o(1)),\qquad
 t(\tau)C_\omega(\tau)\longrightarrow\frac1{\bar c_{\omega,a}}.
\]
Uniform convergence of \(\Omega\) to \(\bar\Omega_a\) and the
reconstruction identity therefore give the decay asymptotic
\begin{equation}\label{eq:physical-decay}
\begin{aligned}
 t(\tau)\|\omega(\cdot,t(\tau))\|_{L^\infty(\bbR)}
 &=t(\tau)C_\omega(\tau)\|\Omega(\cdot,\tau)\|_{L^\infty(\bbR)}\\
 &\longrightarrow
 \frac{\|\bar\Omega_a\|_{L^\infty(\bbR)}}{\bar c_{\omega,a}}
 \qquad(\tau\to\infty).
\end{aligned}
\end{equation}
For a zero amplitude rate, \eqref{eq:scale-asymptotics} gives
\[
 \frac{t(\tau)}\tau\longrightarrow e^{I_\omega},\qquad
 C_\omega(\tau)\longrightarrow e^{-I_\omega}.
\]
For constant rates these formulas reduce exactly to
\eqref{eq:exact-solution} and its exponential limit.

It remains to identify the reconstruction with the classical solution
issued from the smooth data in \eqref{eq:basin}. Indeed, normalize
that solution using the outer characteristic
\(R_{\rm out}'(t)=a u(R_{\rm out}(t),t)\), \(R_{\rm out}(0)=1\),
and set
\(C_l=R_{\rm out}^{-1}\),
\(C_\omega=-\omega_x(0,t)/C_l>0\), and
\(\idiff\tau/\idiff t=C_\omega\).
The slope equation preserves this positivity, and direct
differentiation gives \eqref{eq:normalization}.
The rescaled classical solution belongs to the affine graph class,
so uniqueness identifies it with the global strong solution.
Any finite physical breakdown time before \(t(\infty)\) would
correspond to a finite rescaled time, making
\(C_\omega\|\Omega\|_\infty\) bounded.
Lem.~\ref{lem:real-local} below, including persistence of smoothness,
excludes that breakdown. For \(a<3/2\), the smooth solution therefore
develops a singularity at \(T_{\rm phys}\). For \(a>3/2\), the relation
\(t(\tau)\to\infty\) excludes every finite breakdown time and gives
global smoothness with the decay asymptotic \eqref{eq:physical-decay}. This
completes Thm.~\ref{thm:profiles}.

\section{Signed-data global regularity}\label{sec:global}

The global argument has two independent parts. A positive quadratic
force makes an inverse power of the flow Jacobian convex. This
excludes collapse uniformly over all active labels once the origin
strain is integrable. A common mass and high-moment estimate supplies
that integrability in the three parameter regimes.

\subsection{Real-line continuation and the flow}

The final global step will bound the vorticity amplitude. We first
prove that this bound prevents a finite smooth lifespan, keeping the
origin-normalized velocity outside any unnecessary \(L^2\) assumption.
The characteristic identities then express that amplitude through
an inverse power of the flow Jacobian.

\begin{lemma}[Real-line local theory and continuation]\label{lem:real-local}
For \(a>0\) and odd compactly supported
\(\omega_0\in H^1(\bbR)\), equation \eqref{eq:model} has a unique
maximal solution
\[
 \omega\in C([0,T_*);H^1(\bbR))
       \cap C^1([0,T_*);L^2(\bbR)).
\]
Its support is compact on every compact time interval. If \(T_*\)
is finite and \(\sup_{t<T_*}\|\omega(t)\|_\infty\) is finite,
the solution continues beyond \(T_*\). Every initially finite
higher Sobolev norm persists while the \(H^1\) norm stays bounded.
\end{lemma}
\begin{proof}
Hilbert-transform boundedness on \(H^1\) gives a Lipschitz
origin-normalized velocity with
\(|u(x,t)|\leq |x|\|\bmH_{\bbR}\omega(t)\|_\infty\).
Construct the solution by frozen linear transport iteration with
common initial data,
\[
 \partial_t\omega^{(n+1)}+a u^{(n)}\partial_x\omega^{(n+1)}
 =(\bmH_{\bbR}\omega^{(n)})\omega^{(n+1)},\qquad
 u_x^{(n)}=\bmH_{\bbR}\omega^{(n)},\quad u^{(n)}(0,t)=0.
\]
Smooth odd approximations of the data may be used first.
Differentiating this linear equation and using \(H^1\hookrightarrow
L^\infty\) and the \(L^2\) Hilbert isometry yields
\[
 \frac{\idiff}{\idiff t}\|\omega^{(n+1)}\|_{H^1}
 \leq C_{\rm lin}(a)\|\omega^{(n)}\|_{H^1}
                     \|\omega^{(n+1)}\|_{H^1},
\]
where \(C_{\rm lin}(a)\) is fixed and depends only on \(a\).
This gives a common short-time \(H^1\) bound and a common support
radius \(R_{\rm sup}\), because the characteristics obey the
displayed linear-growth bound.

For completeness, two solutions \(g_1,g_2\) supported in
\([-R_{\rm sup},R_{\rm sup}]\) have vorticity difference
\(g_\Delta=g_1-g_2\) and velocity difference
\(u_\Delta(x)=\int_0^x(\bmH_{\bbR}g_\Delta)(s)\idiff s\).
The inequality
\(\|u_\Delta\|_{L^\infty(-R_{\rm sup},R_{\rm sup})}
\leq\sqrt{R_{\rm sup}}\|g_\Delta\|_2\)
gives the difference estimate
\begin{equation}\label{eq:real-difference}
 \frac12\frac{\idiff}{\idiff t}\|g_\Delta\|_2^2
 \leq\left((1+a/2)\|\bmH_{\bbR}g_1\|_\infty
       +\|g_2\|_\infty
       +a\sqrt{R_{\rm sup}}\|(g_2)_x\|_2\right)\|g_\Delta\|_2^2.
\end{equation}
Indeed, the difference equation has transport velocity from \(g_1\)
and source
\((\bmH_{\bbR}g_1)g_\Delta+
(\bmH_{\bbR}g_\Delta)g_2-a u_\Delta(g_2)_x\).
For consecutive frozen iterates the last two terms contain the
previous coefficient difference. The corresponding two-term
energy bound is a contraction after shortening the interval.
Thus the iterates converge in \(C_tL^2\), with a bounded weak
\(H^1\) limit. The same estimate proves uniqueness.

Interpolation gives continuity into \(H^{3/4}\), and hence
continuous Lipschitz velocity coefficients.
The distributional derivative solves
\begin{equation}\label{eq:derivative-equation}
 (\omega_x)_t+a u(\omega_x)_x
 =(1-a)u_x\omega_x+\omega\bmH_{\bbR}\omega_x.
\end{equation}
Its source is strongly measurable in \(L^\infty_tL^2_x\).
The Lipschitz characteristic formula, with Jacobian bounded
above and below, therefore yields strong \(L^2\) time continuity
of \(\omega_x\), by the same composition and Bochner-integral
argument as for \eqref{eq:weighted-transport}. Hence the solution
has the stated \(C_tH^1\) and \(C_t^1L^2\) regularity.
This construction never requires \(u\in L^2(\bbR)\).

Integration by parts in the equations for \(\omega\) and
\(\omega_x\) gives the energy bound
\begin{equation}\label{eq:real-energy}
 \frac{\idiff}{\idiff t}\|\omega\|_{H^1}
 \leq\left((1+a/2)\|\bmH_{\bbR}\omega\|_\infty
                 +\|\omega\|_\infty\right)\|\omega\|_{H^1}.
\end{equation}
The explicit logarithmic inequality
\begin{equation}\label{eq:log-Hilbert}
 \|\bmH_{\bbR}g\|_\infty
 \leq2+\frac8\pi\|g\|_\infty\log(e+\|g\|_{H^1})
\end{equation}
holds for every \(g\in H^1(\bbR)\). To prove it, split the
Hilbert kernel in its displacement variable \(z\) at the radii
\((e+\|g\|_{H^1})^{-2}\) and \((e+\|g\|_{H^1})^2\).
Subtract \(g(x)\) in the inner principal value.
Cauchy--Schwarz bounds that unnormalized integral by \(4\);
the middle integral is bounded by
\(8\|g\|_\infty\log(e+\|g\|_{H^1})\), and the outer integral
by \(\sqrt2\). Division by \(\pi\) proves \eqref{eq:log-Hilbert}.
If the vorticity supremum is bounded on a finite interval,
\eqref{eq:real-energy} and \eqref{eq:log-Hilbert} give a linear
Gronwall inequality for \(\log(e+\|\omega\|_{H^1})\).
The \(H^1\) norm cannot diverge there. The characteristic formula
for \eqref{eq:derivative-equation} gives a strong \(H^1\) limit
at the endpoint, and the local construction restarts it.

For each integer \(\ell\geq2\), the differentiated product and
transport estimates give
\begin{equation}\label{eq:real-smooth}
 \frac{\idiff}{\idiff t}\|\omega\|_{H^\ell}
 \leq C_{\rm Sob}(a,\ell)\|\omega\|_{H^1}\|\omega\|_{H^\ell},
\end{equation}
with one fixed constant \(C_{\rm Sob}(a,\ell)\) for each order.
Here is the derivative count in the transport commutator. Its
only potentially problematic endpoint product is
\(\|\partial_x^\ell u\|_2\|\omega_x\|_\infty
=\|\partial_x^{\ell-1}\omega\|_2\|\omega_x\|_\infty\).
Interpolation between \(H^1\) and \(H^\ell\), together with
the bound \(\|g\|_\infty^2\leq2\|g\|_2\|g_x\|_2\), gives
\[
 \|\partial_x^\ell u\|_2\|\omega_x\|_\infty
 \leq\sqrt2\|\omega\|_{H^1}^{1+1/(2(\ell-1))}
\|\omega\|_{H^\ell}^{1-1/(2(\ell-1))}
\leq\sqrt2\|\omega\|_{H^1}\|\omega\|_{H^\ell}.
\]
The intermediate Leibniz terms obey the same bound by Hölder and
the same interpolation. The stretching product has the usual
Leibniz bound with one factor in \(L^\infty\); both that factor
and its Hilbert transform are controlled by \(H^1\).
This proves \eqref{eq:real-smooth} and persistence of smoothness.
\end{proof}

With the continuation criterion established, we now specialize to the signed-data class and introduce the Lagrangian variables used in the global argument.
Put \(R_0=R_{\omega_0}\) for nontrivial data and \(R_0=0\)
for zero data. Except when treating that immediate case in
Lem.~\ref{lem:noncollapse}, take a nontrivial datum in \(\bbG\)
for the rest of the global argument. On the positive half-line define
the nonnegative function $v$ by
\begin{equation}\label{eq:positive-vorticity}
 v(x,t)=-\omega(x,t),\qquad x>0.
\end{equation}
Oddness and sign are preserved by uniqueness and characteristics.
Pairing the odd Hilbert kernel gives
\begin{equation}\label{eq:positive-velocity}
 u_x(x,t)=\frac2\pi\operatorname{p.v.}\int_0^\infty
 \frac{s v(s,t)}{s^2-x^2}\idiff s,\qquad
 u(x,t)=\frac1\pi\int_0^\infty v(s,t)
           \log\left|\frac{x+s}{x-s}\right|\idiff s.
\end{equation}
In particular the velocity is nonnegative for \(x>0\).
Define its flow, Jacobian, and inverse-Jacobian power by
\begin{equation}\label{eq:flow}
 \eta_t(t,\alpha)=a u(\eta(t,\alpha),t),\quad
 \eta(0,\alpha)=\alpha,\qquad
 J=\eta_\alpha,\qquad \rho=J^{-1/(2a)}.
\end{equation}
Here \(\alpha\) is a Lagrangian label, and the active labels
lie in \([0,R_0]\) and their odd reflection. The identities
\begin{equation}\label{eq:flow-identities}
 J_t=a u_x(\eta,t)J,\qquad
 \rho_t=-\frac12u_x(\eta,t)\rho,\qquad
 \omega(\eta(t,\alpha),t)=\omega_0(\alpha)J^{1/a}
                        =\frac{\omega_0(\alpha)}{\rho(t,\alpha)^2}
\end{equation}
follow directly along each characteristic.

\subsection{The common Mellin estimate}

The positive-half-line sign turns the origin strain into a positive
moment. We pair that moment with the mass and a higher moment whose
derivative has a favorable sign. Symmetrization and logarithmic
Fourier analysis supply the common estimate; the admissible exponent
range will determine the separate endpoint arguments.

To reveal the sign structure of the moment evolution, we introduce two Mellin kernels and their associated quadratic forms.
Define the moment family and the origin strain by
\begin{equation}\label{eq:moments}
 \frQ_\lambda(t)=\int_0^\infty v(x,t)x^{-\lambda}\idiff x,\qquad
 h(t)=u_x(0,t)=\frac2\pi\frQ_1(t).
\end{equation}
The index \(\lambda\) is a real moment exponent. Ordinary smooth
odd data make these moments finite for \(0\leq\lambda<2\).
A compact-time cubic bound at zero makes them finite for
\(0\leq\lambda<4\). In particular \(\frQ_0\) is the mass.

For \(s>0\) and real kernel order \(m>0\), define
\begin{equation}\label{eq:mellin-kernels}
 \ell_{\rm vel}(s)=\log\left|\frac{s+1}{s-1}\right|,\qquad
 \Psi_m(s)=\frac{2s(s^m-1)}{s^2-1},\qquad
 \Phi_m(s)=(s^m+1)\ell_{\rm vel}(s)-\Psi_m(s).
\end{equation}
The removable value of \(\Psi_m\) at one is \(m\); the
logarithmic singularity of \(\Phi_m\) is retained in its integrals.
The corresponding quadratic forms are
\begin{equation}\label{eq:mellin-forms}
 \begin{split}
 \frA_{\Phi_m}(g)&=\int_0^\infty\int_0^\infty
     g(x)g(s)x^{-m}\Phi_m(x/s)\idiff x\idiff s,\\
 \frA_{\Psi_m}(g)&=\int_0^\infty\int_0^\infty
     g(x)g(s)x^{-m}\Psi_m(x/s)\idiff x\idiff s .
 \end{split}
\end{equation}

\begin{lemma}[Moment identity]\label{lem:moment-identity}
At every regular time the identity
\begin{equation}\label{eq:moment-identity}
 \frQ_\lambda'
 =-\frac{a\lambda}{2\pi}\frA_{\Phi_{\lambda+1}}(v)
  -\frac{a\lambda-a-1}{2\pi}\frA_{\Psi_{\lambda+1}}(v)
\end{equation}
holds for exponents \(0<\lambda<2\). Under the cubic bound
\(v(x,t)=O(x^3)\), uniform on compact regular time intervals,
it also holds for
\(0<\lambda\leq3\).
\end{lemma}
\begin{proof}
The equation \(v_t+a u v_x=u_xv\) and integration by parts give
\[
 \frQ_\lambda'=(a+1)\int_0^\infty u_xv x^{-\lambda}\idiff x
              -a\lambda\int_0^\infty uv x^{-\lambda-1}\idiff x.
\]
The lower cutoff flux is \(O(x^{2-\lambda})\) for ordinary
odd data and \(O(x^{4-\lambda})\) for cubic data. Compact
support removes the other boundary. With \(m=\lambda+1\), first
delete the diagonal symmetrically in the double integral from
\eqref{eq:positive-velocity} and then exchange \(x\) and \(s\).
Averaging the two expressions gives
\begin{equation}\label{eq:moment-symmetrization}
 \begin{split}
 \int_0^\infty u_xv x^{-\lambda}\idiff x
 &=\frac1\pi\int_0^\infty\int_0^\infty
 v(x)v(s)\frac{s x^{-\lambda}-x s^{-\lambda}}{s^2-x^2}
       \idiff x\idiff s
 =\frac1{2\pi}\frA_{\Psi_m}(v),\\
 \int_0^\infty uv x^{-m}\idiff x
 &=\frac1{2\pi}\bigl(\frA_{\Phi_m}(v)+\frA_{\Psi_m}(v)\bigr).
 \end{split}
\end{equation}
These formulas imply \eqref{eq:moment-identity}.

The symmetrized strain kernel has the removable diagonal value
\[
 \lim_{s\to x}
 \frac{s x^{-\lambda}-x s^{-\lambda}}{s^2-x^2}
 =\frac{\lambda+1}{2}x^{-\lambda-1}\qquad(x>0).
\]
Thus the principal-value pole cancels; the velocity kernel retains
only a locally integrable logarithm. To justify the remaining
cutoff limits, fix a compact regular time interval \([0,T']\)
and the exponent \(\lambda\). Smoothness and bounded support give
uniform bounds on \(v(x,t)/x\), or on \(v(x,t)/x^3\) in the
cubic case. On the region \(0<x<s/2\), each absolute double
integrand in \eqref{eq:moment-symmetrization} is consequently
bounded by a fixed multiple of
\[
 x^{1-\lambda}\quad\hbox{for ordinary data},\qquad
 x^{3-\lambda}s^2\quad\hbox{for cubic data}.
\]
On \(s/2<x<s\), the corresponding majorants are
\[
 s^{1-\lambda}\bigl(1+|\log(1-x/s)|\bigr),\qquad
 s^{5-\lambda}\bigl(1+|\log(1-x/s)|\bigr).
\]
The multiples depend on the fixed exponent, data, and regular time
interval. All four majorants are integrable in their asserted
ranges, and symmetry covers \(x>s\). At the cubic endpoint
\(\lambda=3\), the explicit identities
\(\Psi_4(s)=2(s+s^3)\) and
\(\Phi_4(s)=-4s^3/3+O(s^5)\) also exhibit the cancellation.
Dominated convergence now justifies the symmetrization and the
moment identity, including that endpoint. No uniform cutoff bound
as \(\lambda\) approaches two is required.
\end{proof}

The first form in the moment identity need not have a pointwise
positive kernel. The next lemma proves positivity of the quadratic
form through its logarithmic Fourier symbol. A separate secant bound
for the second kernel then converts it into a product of moments.

\begin{lemma}[Mellin positivity and the secant bound]\label{lem:mellin-positive}
For kernel orders \(2<m\leq4\), the form \(\frA_{\Phi_m}\) is
nonnegative on its logarithmic \(L^2\) domain, and the inequality
\begin{equation}\label{eq:secant}
 \Psi_m(s)\geq\frac m2(s+s^{m-1})\qquad(s>0)
\end{equation}
holds.
The form domain consists of real inputs for which
\(G_m(g)(r)=g(e^r)e^{(1-m/2)r}\) belongs to \(L^2(\bbR)\).
\end{lemma}
\begin{proof}
Under logarithmic substitution, define the convolution kernel by
\begin{equation}\label{eq:mellin-convolution}
 W_m(z)=e^{-mz/2}\Phi_m(e^z)
 =2\cosh(mz/2)\log|\coth(z/2)|
       -2\frac{\sinh(mz/2)}{\sinh z}.
\end{equation}
The quadratic form becomes
\(\int_{\bbR}\int_{\bbR}G_m(g)(r)G_m(g)(s)W_m(r-s)\idiff r\idiff s\).
With the roman symbol \(\iunit\) denoting the imaginary unit,
we use the Fourier convention
\begin{equation}\label{eq:fourier}
 \widehat g(\xi)=\int_{\bbR}g(z)e^{-\iunit\xi z}\idiff z,\qquad
 g(z)=\frac1{2\pi}\int_{\bbR}\widehat g(\xi)e^{\iunit\xi z}\idiff\xi.
\end{equation}
The same Fourier notation denotes the extension to tempered
distributions when an input is not integrable.
The kernel symbol is
\begin{equation}\label{eq:mellin-symbol}
 \widehat W_m(\xi)=
 \frac{2\pi\left(\xi\sinh(\pi\xi)
  -(\xi^2+m^2/4-m/2)\sin(\pi m/2)\right)}
 {(\xi^2+m^2/4)(\cosh(\pi\xi)+\cos(\pi m/2))}.
\end{equation}
Here is a derivation without a transform table. Let \(w\) denote
a complex argument distinct from the real integration variables.
For a positive argument \(z\),
the expansion
\(\log|\coth(z/2)|=2\sum_{j\geq1,\ j\ {\rm odd}}e^{-jz}/j\)
and the partial fraction identity
\[
 \sum_{\substack{j\geq1\\j\ {\rm odd}}}
 \left(\frac1{j-w}-\frac1{j+w}\right)
 =\frac\pi2\tan\frac{\pi w}{2}
\]
give
\(\big(\log|\coth(z/2)|\big)^{\wedge}(\xi)
=4\sum_{j\ {\rm odd}}(j^2+\xi^2)^{-1}
=\pi\tanh(\pi\xi/2)/\xi\).
The partial fraction identity follows by logarithmic differentiation
of the locally uniformly convergent product
\(\cos(\pi w/2)=\prod_{j\geq1,\ j\ {\rm odd}}(1-w^2/j^2)\);
it holds away from the odd integers.
For \(0<m<2\), expansion of \(1/\sinh z\) gives in the same way
\[
 {\big(\sinh(mz/2)/\sinh z\big)^{\wedge}}(\xi)
 =\frac{\pi\sin(\pi m/2)}
        {\cosh(\pi\xi)+\cos(\pi m/2)}.
\]
For \(0<m<2\), the complex frequency shifts remain within the
strip of convergence of the transform of \(\log|\coth(z/2)|\).
They give the explicit formula
\[
 \begin{aligned}
 &\big(2\cosh(mz/2)\log|\coth(z/2)|\big)^{\wedge}(\xi)\\
 &\quad=\pi\left(
 \frac{\tanh(\pi(\xi+\iunit m/2)/2)}{\xi+\iunit m/2}
 +\frac{\tanh(\pi(\xi-\iunit m/2)/2)}{\xi-\iunit m/2}\right)\\
 &\quad=\frac{2\pi\bigl(\xi\sinh(\pi\xi)
                   +(m/2)\sin(\pi m/2)\bigr)}
 {\bigl(\xi^2+m^2/4\bigr)
  \bigl(\cosh(\pi\xi)+\cos(\pi m/2)\bigr)}.
 \end{aligned}
\]
The last equality uses
\(\tanh(A+\iunit B)=(\sinh(2A)+\iunit\sin(2B))
/(\cosh(2A)+\cos(2B))\) for real \(A,B\) wherever the
denominator is nonzero. Subtracting twice the preceding transform
of \(\sinh(mz/2)/\sinh z\) now gives \eqref{eq:mellin-symbol}.

To justify continuation past \(m=2\) and include \(m=4\), use
the combined series
\begin{equation}\label{eq:mellin-series}
 W_m(z)=2\sum_{\substack{j\geq1\\j\ {\rm odd}}}
 \left((j^{-1}-1)e^{-(j-m/2)z}
             +(j^{-1}+1)e^{-(j+m/2)z}\right),\qquad z>0.
\end{equation}
The first coefficient at \(j=1\) vanishes. Near zero the kernels
are bounded on compact \(m\)-ranges by a fixed multiple of
\(1+|\log z|\), and at infinity the first surviving shifted
exponential starts at \(j=3\). Thus they form an \(L^1\)-analytic
family on \(0<\operatorname{Re}m<6\), which continues the symbol
identity and permits its limit at four. In the asserted range
\(2<m\leq4\), the factor \(\sin(\pi m/2)\) is nonpositive,
and \(m^2/4-m/2\) is nonnegative. The symbol is therefore
nonnegative, including the explicit endpoint
\begin{equation}\label{eq:mellin-four}
 \widehat W_4(\xi)=
 \frac{2\pi\xi\sinh(\pi\xi)}
      {(\xi^2+4)(\cosh(\pi\xi)+1)}\geq0.
\end{equation}
Plancherel's identity therefore gives the positivity formula
\[
 \frA_{\Phi_m}(g)
 =\frac1{2\pi}\int_{\bbR}\widehat W_m(\xi)
      |\widehat{G_m(g)}(\xi)|^2\idiff\xi\geq0.
\]

Finally, the fundamental theorem of calculus gives
\[
 \Psi_m(s)=ms\int_0^1(1+r(s^2-1))^{m/2-1}\idiff r.
\]
Since \(0<m/2-1\leq1\), concavity gives the pointwise lower bound
\[
 (1+r(s^2-1))^{m/2-1}
 \geq(1-r)+r s^{m-2}\qquad(0\leq r\leq1).
\]
Integrating this inequality yields
\(\Psi_m(s)\geq ms(1+s^{m-2})/2\), which is
\eqref{eq:secant}. All applications satisfy
the form-domain hypothesis: ordinary data give
\(G_m(v)(r)=O(e^{(2-m/2)r})\) with \(m<3\) at negative
infinity, and cubic data give \(G_4(v)(r)=O(e^{2r})\).
Compact support controls positive infinity.
\end{proof}

The high moment decreases under the coefficient condition below,
whereas the mass may increase. Their interpolation product controls
the origin strain and combines these two differential bounds into
one estimate usable in every global regime.

\begin{proposition}[One mass and high-moment estimate]\label{prop:moment}
Fix a nontrivial regular signed solution and a moment exponent
\(\lambda\). Assume either the range \(1<\lambda<2\), or the
range \(1<\lambda\leq3\) together with compact-time cubic vanishing
at the origin. Under the condition \(a\lambda-a-1\geq0\),
its energy
\begin{equation}\label{eq:moment-energy}
 \frE_\lambda=\frQ_\lambda^{1/\lambda}
                     \frQ_0^{1-1/\lambda}
\end{equation}
satisfies
\begin{equation}\label{eq:common-moment}
 \frac{\frQ_\lambda'}{\frQ_\lambda}
 \leq-\frac{(\lambda+1)(a\lambda-a-1)}4h,\qquad
 \frac{\frQ_0'}{\frQ_0}\leq\frac{a+1}4h,\qquad
 h\leq\frac2\pi\frE_\lambda,\qquad
 \frac{\frE_\lambda'}{\frE_\lambda}
 \leq\frac{a+2-a\lambda}{4}h .
\end{equation}
\end{proposition}
\begin{proof}
Lemma~\ref{lem:mellin-positive} makes the first term in
\eqref{eq:moment-identity} nonpositive. The secant bound
\eqref{eq:secant} and the sign of \(v\) give
\[
 \frA_{\Psi_{\lambda+1}}(v)
 \geq(\lambda+1)\frQ_\lambda\frQ_1
 =\frac{\pi(\lambda+1)}2h\frQ_\lambda,
\]
proving the first inequality in \eqref{eq:common-moment}.
For the mass, direct symmetrization gives
\[
 \frQ_0'=(a+1)\int_0^\infty u_xv\idiff x
 =\frac{a+1}{\pi}\int_0^\infty\int_0^\infty
                    \frac{v(x)v(s)}{x+s}\idiff x\idiff s .
\]
The elementary inequality
\((x+s)^{-1}\leq(x^{-1}+s^{-1})/4\) proves the mass bound.
Hölder gives \(\frQ_1\leq
\frQ_\lambda^{1/\lambda}\frQ_0^{1-1/\lambda}\).
The logarithmic derivative of \eqref{eq:moment-energy} consequently
satisfies
\[
 \begin{aligned}
 \frac{\frE_\lambda'}{\frE_\lambda}
 &=\frac1\lambda\frac{\frQ_\lambda'}{\frQ_\lambda}
   +\left(1-\frac1\lambda\right)\frac{\frQ_0'}{\frQ_0}\\
 &\leq\frac{-(\lambda+1)(a\lambda-a-1)
                  +(\lambda-1)(a+1)}{4\lambda}\,h
 =\frac{a+2-a\lambda}{4}\,h,
 \end{aligned}
\]
which is the last bound in \eqref{eq:common-moment}.
All moments are positive for nontrivial signed data, so these
logarithmic derivatives are legitimate.
\end{proof}

\subsection{A positive force and uniform noncollapse}

The moment bound concerns the origin alone. To control vorticity
throughout the support, we derive a nonnegative quadratic force in
the material strain equation. It makes the inverse-Jacobian power
convex and supplies a fixed-label barrier. The final compactness
argument upgrades that barrier to a bound uniform over all labels.

For a real input \(g\) whose Fourier transform satisfies
\(\int_{\bbR}(1+\xi^2)|\widehat g(\xi)|\idiff\xi<\infty\),
define the pointwise quadratic force by
\begin{equation}\label{eq:force}
 \frF(g)=-g g_{xx}
         -\bmH_{\bbR}\bigl(g\bmH_{\bbR}g_{xx}\bigr).
\end{equation}
It has the nonnegative square representation
\begin{equation}\label{eq:force-square}
 \frF(g)(x)=\frac1{\pi^2}\int_0^\infty
 \xi\left|\int_\xi^\infty\widehat g(\nu)e^{\iunit\nu x}\idiff\nu\right|^2
 \idiff\xi\geq0.
\end{equation}
To prove it, the Hilbert multiplier \(-\iunit\operatorname{sgn}\xi\)
gives the double Fourier kernel
\(\nu^2(1-\operatorname{sgn}(\xi+\nu)\operatorname{sgn}\nu)\)
with factor \(1/(4\pi^2)\).
Splitting the two opposite-sign regions and using reality
reduces the expression to
\[
 \frac1{2\pi^2}\int_0^\infty\int_0^\infty
 \min(r,s)^2\widehat g(r)\widehat g(-s)e^{\iunit(r-s)x}\idiff r\idiff s.
\]
The identity \(\min(r,s)^2=2\int_0^{\min(r,s)}z\idiff z\)
and Fubini prove \eqref{eq:force-square}.
Absolute convergence follows from the assumed weighted Fourier
\(L^1\) norm, so smooth Fourier approximation passes this
calculation to the stated class.

For a smooth compactly supported odd vorticity, the identity
\(\widehat u(\xi)=-\widehat\omega(\xi)/|\xi|\) holds at nonzero
frequencies. The numerator vanishes linearly at zero, and
its far-field decay is rapid. Thus the force hypothesis holds
for the physical velocity. Its inverse Fourier transform is
odd and has \(u(0)=0\), agreeing with the chosen normalization.

We next derive the material strain identity in which this nonnegative force enters the dynamics. The Cotlar identity
\[
 \bmH_{\bbR}(\omega\bmH_{\bbR}\omega)
 =\frac12\bigl((\bmH_{\bbR}\omega)^2-\omega^2\bigr)
\]
follows by squaring the positive-frequency function
\(\omega+\iunit\bmH_{\bbR}\omega\) and comparing its real and
imaginary parts. Since
\(\bmH_{\bbR}u_{xx}=-\omega_x\), the force definition gives
\(\frF(u)=-u u_{xx}+\bmH_{\bbR}(u\omega_x)\).
Applying \(\bmH_{\bbR}\) to \eqref{eq:model} and adding the
material transport term therefore gives
\begin{equation}\label{eq:strain-material}
 \begin{aligned}
 (\partial_t+a u\partial_x)u_x
 &=\bmH_{\bbR}(\omega\bmH_{\bbR}\omega)
     -a\bmH_{\bbR}(u\omega_x)+a u u_{xx}\\
 &=\frac12(u_x^2-\omega^2)-a\frF(u).
 \end{aligned}
\end{equation}
In differentiating \(\rho_t=-u_x(\eta,t)\rho/2\), this identity
gives the cancellation
\[
 \begin{aligned}
 \rho_{tt}
 &=-\frac12\bigl((\partial_t+a u\partial_x)u_x\bigr)(\eta,t)\rho
     +\frac14u_x(\eta,t)^2\rho\\
 &=\frac14\omega(\eta,t)^2\rho
     +\frac a2\frF(u)(\eta)\rho.
 \end{aligned}
\]
The two terms containing \(u_x^2\) have canceled. Substituting
the vorticity identity in \eqref{eq:flow-identities} then yields
\begin{equation}\label{eq:convexity}
 \rho_{tt}(t,\alpha)
 =\frac{\omega_0(\alpha)^2}{4\rho(t,\alpha)^3}
   +\frac a2\frF(u)(\eta(t,\alpha))\rho(t,\alpha)\geq0.
\end{equation}
The sign in this identity is the reason a pointwise force estimate
is useful. Positivity at each fixed label is still insufficient near
labels with zero initial vorticity. The next lemma treats those
labels separately and uses the common convexity bound to prove
uniformity without assuming regularity at the terminal time.

\begin{lemma}[Uniform noncollapse]\label{lem:noncollapse}
Let the parameter \(a\) be positive, let the initial datum
\(\omega_0\) belong to \(\bbG\), and let a smooth solution
exist on a finite interval \([0,T)\). Under the integrability
condition \(\int_0^T h(t)\idiff t<\infty\), the bounds
\begin{equation}\label{eq:noncollapse}
 \inf_{\substack{0\leq t<T\\0\leq\alpha\leq R_0}}
       \rho(t,\alpha)>0,\qquad
 \sup_{0\leq t<T}\|\omega(t)\|_\infty<\infty
\end{equation}
hold.
\end{lemma}
\begin{proof}
For zero data the convention \(R_0=0\) and the identity
\(\rho=1\) give both conclusions. Assume henceforth that the
data are nontrivial, and write \(u_0\) for their normalized velocity.
The initial values
\[
 \rho(0,\alpha)=1,\qquad
 \rho_t(0,\alpha)=-\frac12u_{0x}(\alpha)
\]
follow from \eqref{eq:flow} and \eqref{eq:flow-identities}.
We first obtain a positive lower bound at each fixed label.
We then use convexity to make those bounds uniform.

Fix a label \(\alpha\) with \(\omega_0(\alpha)\ne0\).
Its barrier energy is defined by
\begin{equation}\label{eq:barrier}
 \frB_\alpha^{\mathrm{Lag}}(t)
 =\rho_t(t,\alpha)^2+
       \frac{\omega_0(\alpha)^2}{4\rho(t,\alpha)^2}.
\end{equation}
All denominators are positive on the regular interval.
Keeping the label fixed, differentiation and
\eqref{eq:convexity} give the cancellation
\[
 \begin{aligned}
 \frac{\idiff}{\idiff t}\frB_\alpha^{\mathrm{Lag}}
 &=2\rho_t\rho_{tt}
       -\frac{\omega_0(\alpha)^2}{2\rho^3}\rho_t\\
 &=2\rho_t\left(
       \frac{\omega_0(\alpha)^2}{4\rho^3}
       +\frac a2\frF(u)(\eta(t,\alpha))\rho\right)
       -\frac{\omega_0(\alpha)^2}{2\rho^3}\rho_t\\
 &=a\frF(u)(\eta(t,\alpha))\rho\rho_t .
 \end{aligned}
\]
Here the unlabelled occurrences of \(\rho\) and its derivatives
are evaluated at \((t,\alpha)\).
The nonnegativity of the force therefore makes the energy
nonincreasing whenever \(\rho_t\leq0\).

Convexity makes \(\rho_t\) nondecreasing.
If \(\rho_t(0,\alpha)\geq0\), the bound
\(\rho(t,\alpha)\geq1\) holds throughout the interval.
Otherwise, the decreasing interval begins at zero and lasts until
the derivative first reaches zero, if it does so before \(T\).
During that interval, the comparison
\[
 \frac{\omega_0(\alpha)^2}{4\rho(t,\alpha)^2}
 \leq\frB_\alpha^{\mathrm{Lag}}(t)
 \leq\frB_\alpha^{\mathrm{Lag}}(0)
 =\frac{u_{0x}(\alpha)^2+\omega_0(\alpha)^2}{4}
\]
gives the barrier
\begin{equation}\label{eq:pointwise-barrier}
 \rho(t,\alpha)^2\geq
 \frac{\omega_0(\alpha)^2}
      {u_{0x}(\alpha)^2+\omega_0(\alpha)^2}.
\end{equation}
If a minimum occurs before \(T\), subsequent monotonicity preserves
this bound. If the decreasing interval lasts until \(T\), the
comparison already covers every regular time.
The right side of \eqref{eq:pointwise-barrier} is positive and
at most one, so the same bound also holds in the initially
nondecreasing case.

The remaining labels require the hypotheses on the zero set.
At an interior zero \(z\), the compressive-zero condition gives
\(\rho_t(0,z)=-u_{0x}(z)/2\geq0\).
Convexity consequently gives \(\rho(t,z)\geq1\).
At the origin, the fixed characteristic \(\eta(t,0)=0\) gives
the identity
\begin{equation}\label{eq:origin-rho}
 \rho(t,0)=
 \exp\left(-\frac12\int_0^t h(s)\idiff s\right).
\end{equation}
Since \(h\geq0\), the assumed integrability gives the uniform
lower bound
\[
 \rho(t,0)\geq
 \exp\left(-\frac12\int_0^T h(s)\idiff s\right)>0.
\]
At the support endpoint, the initial strain has the representation
\[
 u_{0x}(R_0)=\frac2\pi\int_0^{R_0}
       \frac{s v(s,0)}{s^2-R_0^2}\idiff s<0.
\]
Smooth compact support makes the integrand integrable at \(R_0\);
its sign is nonpositive, and nontriviality makes the integral
strictly negative. Thus \(\rho_t(0,R_0)>0\), and convexity gives
\(\rho(t,R_0)\geq1\).

The fixed-label bounds do not directly give a uniform bound:
the right side of \eqref{eq:pointwise-barrier} can degenerate as
the label approaches a zero of the initial vorticity.
Suppose, for contradiction, that uniform noncollapse fails.
Choose labels \(\alpha_n\in[0,R_0]\) and times \(t_n<T\)
such that \(\rho(t_n,\alpha_n)\to0\).
For each \(n\), choose a minimizing time \(s_n\in[0,t_n]\)
for which the equality
\[
 \rho(s_n,\alpha_n)
   =\min_{0\leq t\leq t_n}\rho(t,\alpha_n)
\]
holds. These minimum values also tend to zero.
After passage to a subsequence, the labels satisfy
\(\alpha_n\to\alpha_{\rm lim}\in[0,R_0]\).

For every fixed \(T'<T\), the smooth flow has a positive,
continuous Jacobian on the compact set
\([0,T']\times[0,R_0]\).
Its inverse power therefore has a positive minimum there.
The minimizing times must consequently satisfy \(s_n\to T\).
Discarding finitely many terms ensures
\(\rho(s_n,\alpha_n)<1\), and hence \(s_n>0\).

Convexity gives the common lower slope bound, while minimality
gives the upper bound in
\[
 -\frac12\|u_{0x}\|_\infty
 \leq\rho_t(0,\alpha_n)
 \leq\rho_t(t,\alpha_n)\leq0
 \qquad(0\leq t\leq s_n).
\]
The last inequality remains valid when \(s_n=t_n\):
a positive derivative at or before \(s_n\), together with
monotonicity of the derivative, would contradict minimality.
Integrating the common slope bound backwards from \(s_n\)
gives the comparison
\[
 \begin{aligned}
 \rho(t,\alpha_n)
 &=\rho(s_n,\alpha_n)
       -\int_t^{s_n}\rho_t(s,\alpha_n)\idiff s\\
 &\leq\rho(s_n,\alpha_n)
       +\frac12\|u_{0x}\|_\infty(s_n-t)
       \qquad(0\leq t\leq s_n).
 \end{aligned}
\]

Fix a time \(t<T\) before passing to the limit.
Eventually \(t<s_n\), and continuity with respect to the label
at this fixed regular time gives the bound
\[
 \rho(t,\alpha_{\rm lim})
 \leq\frac12\|u_{0x}\|_\infty(T-t).
\]
No continuity at \(T\) is used in this passage.
The bound holds for every \(t<T\).
If the limiting label is an interior label with nonzero initial
vorticity, letting \(t\) approach \(T\) contradicts
\eqref{eq:pointwise-barrier}.
An interior zero or the support endpoint instead has the lower
bound one. At the origin, \eqref{eq:origin-rho} has the fixed
positive lower bound obtained above.
These alternatives exhaust \([0,R_0]\), proving the first
conclusion of \eqref{eq:noncollapse}.

The flow preserves zero initial vorticity at every exterior
label, and oddness covers the negative half-line.
The last identity in \eqref{eq:flow-identities} therefore gives
the explicit bound
\[
 \sup_{0\leq t<T}\|\omega(t)\|_\infty
 \leq
 \frac{\|\omega_0\|_\infty}
 {\left(
   \inf_{\substack{0\leq t<T\\0\leq\alpha\leq R_0}}
          \rho(t,\alpha)
  \right)^2}
 <\infty.
\]
This is the continuation quantity required by
Lem.~\ref{lem:real-local}; in particular, \(T\) cannot be a
finite maximal smooth lifespan.
\end{proof}

\subsection{The three moment regimes}

We now supply the origin-strain bound and finish
Thm.~\ref{thm:global}. Suppose that the maximal smooth lifespan
\(T_*\) is finite. A uniform bound for \(h\) on \([0,T_*)\)
will give Lem.~\ref{lem:noncollapse} and then contradict
Lem.~\ref{lem:real-local}. The argument now splits into the case \(a>2\), the endpoint \(a=2\), and the origin-degenerate case \(a\geq1\).

The last coefficient in \eqref{eq:common-moment} is nonpositive
when \(\lambda\geq1+2/a\). For \(a>2\), we can meet this
condition within the admissible range \(1<\lambda<2\).
Choose the midpoint \(\lambda=3/2+1/a\) of \((1+2/a,2)\).
Then \(3/2<\lambda<2\) and \(a\lambda-a-1=a/2>0\), so
Proposition~\ref{prop:moment} gives
\begin{equation}\label{eq:supercritical-moment}
 \frac{\frE_\lambda'}{\frE_\lambda}
 \leq-\frac{a-2}{8}h,\qquad
 0\leq h(t)\leq\frac2\pi\frE_\lambda(0).
\end{equation}
This is the required bound for every finite time interval.

At \(a=2\), the same coefficient condition would require
\(\lambda\geq2\), whereas a general smooth odd datum only
ensures finiteness of \(\frQ_\lambda\) for \(\lambda<2\).
We therefore choose \(\lambda=2-\epsilon\), with
\(0<\epsilon<1/2\). Here \(\epsilon\) is an auxiliary moment
parameter, distinct from the stability threshold
\(\varepsilon_{\rm st}\). The common estimate gives
\begin{equation}\label{eq:critical-moment}
 \frE_{2-\epsilon}'\leq\frac{\epsilon}{\pi}
                         \frE_{2-\epsilon}^2.
\end{equation}
Although the initial energy can diverge in this limit, smooth
odd data satisfy the explicit bound
\[
 \frQ_{2-\epsilon}(0)
 \leq\|\omega_{0x}\|_\infty\frac{R_0^\epsilon}{\epsilon}.
\]
The initial energy therefore satisfies
\begin{equation}\label{eq:critical-initial}
 \epsilon\frE_{2-\epsilon}(0)
 \leq\|\omega_{0x}\|_\infty^{1/(2-\epsilon)}
 R_0^{\epsilon/(2-\epsilon)}
 \frQ_0(0)^{(1-\epsilon)/(2-\epsilon)}
 \epsilon^{(1-\epsilon)/(2-\epsilon)}
 \longrightarrow0.
\end{equation}
The final exponent is at least \(1/3\), and all other factors
stay bounded. We may therefore fix a sufficiently small \(\epsilon\)
for which the bound \(\epsilon\frE_{2-\epsilon}(0)T_*/\pi<1/2\)
holds.
Integrating \eqref{eq:critical-moment} in reciprocal form gives
\[
 \frac1{\frE_{2-\epsilon}(t)}
 \geq\frac1{\frE_{2-\epsilon}(0)}-\frac\epsilon\pi t
 \geq\frac1{2\frE_{2-\epsilon}(0)}
 \qquad(0\leq t<T_*).
\]
Consequently, the bounds
\begin{equation}\label{eq:critical-bound}
 \frE_{2-\epsilon}(t)\leq2\frE_{2-\epsilon}(0),\qquad
 0\leq h(t)\leq\frac4\pi\frE_{2-\epsilon}(0)
\end{equation}
hold on the entire hypothetical lifespan.

For the degenerate assertion, evaluation of
\eqref{eq:derivative-equation} at zero gives
\begin{equation}\label{eq:slope}
 \omega_x(0,t)=\omega_{0x}(0)
       \exp\left((1-a)\int_0^t h(s)\idiff s\right).
\end{equation}
Zero initial slope is preserved. Oddness also gives
\(\omega_{xx}(0,t)=0\), so Taylor's theorem yields, for each
compact regular interval \([0,T']\),
\begin{equation}\label{eq:cubic}
 |\omega(x,t)|\leq\frac{|x|^3}{6}
       \sup_{0\leq s\leq T'}\|\omega_{xxx}(s)\|_\infty
       \qquad(0\leq t\leq T').
\end{equation}
This is precisely the cubic hypothesis in
Lem.~\ref{lem:moment-identity}. Choosing the fixed exponent
\(\lambda=3\) in Prop.~\ref{prop:moment} gives
\begin{equation}\label{eq:degenerate-moment}
 \frac{\frQ_3'}{\frQ_3}\leq-(2a-1)h,\qquad
 \frac{\frE_3'}{\frE_3}\leq-\frac{a-1}{2}h,\qquad
 0\leq h(t)\leq\frac2\pi\frE_3(0)
\end{equation}
for every \(a\geq1\). At \(a=1\) the same fixed energy is
nonincreasing; no exponent limit is used.

The bounds \eqref{eq:supercritical-moment},
\eqref{eq:critical-bound}, and \eqref{eq:degenerate-moment}
show, in their respective regimes, that
\(\int_0^{T_*}h(t)\idiff t<\infty\).
Lemma~\ref{lem:noncollapse} bounds the vorticity
supremum uniformly, including every admissible interior zero.
Lem.~\ref{lem:real-local} extends the solution and preserves
all Sobolev orders, contradicting maximality. This proves both
assertions of Thm.~\ref{thm:global}. The distinct values
\(3/2,2,1\) describe, respectively, the local branch sign and
the two moment endpoints established here.

\vspace{0.2in}
\noindent{\bf Acknowledgments.}
This work was supported in part by NSF Grant DMS-2508463,
the Choi Family Gift Fund, and the Mike Yan Gift Fund.

\appendix
\section{Analytic matrix and row bounds}\label{app:rows}

This appendix completes the matrix estimates used in
Proposition~\ref{prop:linear} and the trace estimate
\eqref{eq:trace-bound}.  The De Gregorio part contributes
\(3/16\) of dissipation in the first weighted mode and at least
\(1/4\) in every other mode, by \eqref{eq:DG-energy}.
To retain the coercivity constant \(3/16\), it therefore suffices
to bound the first symmetric row of the compact remainder by zero
and each remaining row by \(1/16\).
Lemma~\ref{lem:rows} proves precisely these bounds.  Before estimating
the rows, we derive the matrix from the action of \(\calR\) on sine
modes and prove the square summability needed for compactness.

\subsection{Sine coefficients and the weighted-basis matrix}

The coefficient families below distinguish three operations.
The coefficients \(\beta_j^{\rm self}\) and
\(\beta_j^{\rm tail}\) describe the same-frequency and lower-frequency
parts of \(\calR\sin(j\theta)\); \(\beta_j^{\rm tr}\) is its
origin derivative.  Passing from sine modes to
\(e_j=\sin((j+1)\theta)/(j+1)-\sin(j\theta)/j\)
then produces the trace difference
\(\beta_j^{\rm tr\text{-}diff}\) and the adjacent tail sum
\(\beta_j^{\rm tail\text{-}sum}\).  Their exact values are
\begin{equation}\label{eq:beta}
 \begin{gathered}
 \beta_j^{\rm self}=\frac{j-1}{j(2j+1)},\qquad
 \beta_j^{\rm tail}=\frac6{4j^2-1}\quad(j\geq2),\qquad
 \beta_1^{\rm tail}=0,\\
 \beta_j^{\rm tr}=-\frac12-\frac{3(-1)^j}{2(4j^2-1)},\qquad
 \beta_j^{\rm tr\text{-}diff}
 =\frac{\beta_{j+1}^{\rm tr}}{j+1}-\frac{\beta_j^{\rm tr}}j,\qquad
 \beta_j^{\rm tail\text{-}sum}
 =\frac{\beta_{j+1}^{\rm tail}}{j+1}+\frac{\beta_j^{\rm tail}}j.
 \end{gathered}
\end{equation}
We now derive these coefficients and their role in the matrix.
Since \(\bmH_{\bbT}\sin(j\theta)=-\cos(j\theta)\),
\eqref{eq:circle-functionals} gives
\[
 \mathfrak u_1(\sin(j\theta))
 =-\frac14\int_0^\pi\cos(s/2)\cos(js)\idiff s
 =\frac{(-1)^j}{2(4j^2-1)}.
\]
The integrating factor in \eqref{eq:P-ode} therefore yields
\[
 \begin{aligned}
 \cos(\theta/2)\calP\sin(j\theta)
 &=-\frac12\int_0^\theta\cos(s/2)\cos(js)\idiff s
   -\frac{(-1)^j\sin(\theta/2)}{4j^2-1}\\
 &=-\frac{\sin((j-1/2)\theta)}{2(2j-1)}
   -\frac{\sin((j+1/2)\theta)}{2(2j+1)}
   -\frac{(-1)^j\sin(\theta/2)}{4j^2-1}.
 \end{aligned}
\]
Expanding the half-integer sines after division by
\(\cos(\theta/2)\) gives the finite sine expansion
\begin{equation}\label{eq:P-sine}
 \calP\sin(j\theta)=-\frac{\sin(j\theta)}{2j+1}
 +\frac2{4j^2-1}\sum_{i=1}^{j-1}(-1)^{j-i}\sin(i\theta).
\end{equation}
The origin-normalized primitive of
\(\bmH_{\bbT}\sin(j\theta)\) is \(-\sin(j\theta)/j\).
Consequently, \eqref{eq:R} becomes
\begin{equation}\label{eq:R-sine}
 \calR\sin(j\theta)=-\beta_j^{\rm self}\sin(j\theta)
 +\beta_j^{\rm tail}\sum_{i=1}^{j-1}(-1)^{j-i}\sin(i\theta).
\end{equation}
The elementary alternating sum
\[
 \sum_{i=1}^{j-1}i(-1)^{j-i}
 =\begin{cases}
 -j/2,&j\text{ even},\\
 -(j-1)/2,&j\text{ odd}
 \end{cases}
\]
then gives
\[
 (\calR\sin(j\theta))_\theta(0)
 =-j\beta_j^{\rm self}
   +\beta_j^{\rm tail}\sum_{i=1}^{j-1}i(-1)^{j-i}
 =\beta_j^{\rm tr}.
\]
In particular, \eqref{eq:trace-functional} and the definition of
\(e_j\) identify the coefficients used in the transversality estimate:
\begin{equation}\label{eq:trace-basis}
 \mathfrak l(e_j)=\beta_j^{\rm tr\text{-}diff}.
\end{equation}
Applying \eqref{eq:R-sine} to the two sine terms in \(e_j\)
also gives the expansion
\begin{equation}\label{eq:R-basis}
 \begin{aligned}
 \calR e_j={}&-\frac{\beta_{j+1}^{\rm self}}{j+1}
                 \sin((j+1)\theta)
 +\left(\frac{\beta_j^{\rm self}}j
              -\frac{\beta_{j+1}^{\rm tail}}{j+1}\right)\sin(j\theta)\\
 &+\beta_j^{\rm tail\text{-}sum}
       \sum_{i=1}^{j-1}(-1)^{j+1-i}\sin(i\theta).
 \end{aligned}
\end{equation}
This identity explains the adjacent tail sum in \eqref{eq:beta}.

To convert \eqref{eq:R-basis} to the matrix of \(\calK\), let
\(d_k\) denote the ordinary sine coefficients of an odd sine
polynomial with zero origin slope.  If \(b_i\) are its weighted-basis
coefficients, comparison with \eqref{eq:basis} gives
\begin{equation}\label{eq:sine-to-basis}
 k d_k=b_{k-1}-b_k,\qquad b_0=0,\qquad
 b_i=-\sum_{k=1}^i k d_k=\sum_{k=i+1}^\infty k d_k.
\end{equation}
The last equality uses the zero-slope identity \(\sum_k k d_k=0\).
Apply this rule to the zero-slope polynomial
\[
 \calK e_j
 =\frac12\left((\calR e_j)\cos\theta
      -\beta_j^{\rm tr\text{-}diff}\sin\theta\right).
\]
Product-to-sum in \eqref{eq:R-basis} and
\eqref{eq:sine-to-basis} give
\begin{align}
 K_{ij}&=\frac1{4j(j+1)}
 +\frac{(-1)^{i+j+1}(2i+1)}8\beta_j^{\rm tail\text{-}sum}
 &&(i\leq j-2),\label{eq:matrix-dense}\\
 K_{j-1,j}&=\frac14\beta_j^{\rm tail}
 -\frac{\beta_{j+1}^{\rm tail}}{4(j+1)}
 +\frac{j+1}{4j}\beta_j^{\rm self}
 -\frac12\beta_{j+1}^{\rm self}&&(j\geq2),\label{eq:matrix-upper}\\
 K_{jj}&=-\frac{4j^4+8j^3+10j^2+8j+3}
 {4j^2(j+1)^2(2j+1)(2j+3)},\qquad
 K_{j+1,j}=-\frac{j(j+2)}{4(j+1)^2(2j+3)},
 \label{eq:matrix-diagonal}\\
 K_{ij}&=0&&(i\geq j+2).
 \label{eq:matrix-zero}
\end{align}
For example, the highest sine frequency in \(\calK e_j\) is
\(j+2\), with coefficient
\(-\beta_{j+1}^{\rm self}/(4(j+1))\).  Thus
\eqref{eq:sine-to-basis} gives
\[
 K_{j+1,j}
 =-\frac{j+2}{4(j+1)}\beta_{j+1}^{\rm self}
 =-\frac{j(j+2)}{4(j+1)^2(2j+3)}.
\]
For the dense part, the same partial-sum rule gives the useful
alternative formula
\begin{equation}\label{eq:matrix-dense-trace}
 K_{ij}=\frac12\beta_j^{\rm tr\text{-}diff}
 +\frac{(-1)^{j+1}}4\beta_j^{\rm tail\text{-}sum}
 \begin{cases}
 i+1,&i\text{ even},\\
 -i,&i\text{ odd},
 \end{cases}
 \qquad i\leq j-2.
\end{equation}
The trace-difference identity in \eqref{eq:beta-bounds} below reduces
\eqref{eq:matrix-dense-trace} to \eqref{eq:matrix-dense}.

\subsection{Coefficient bounds and compactness}

We require two bounds on these exact coefficients.  The trace bound
controls \(\mathfrak l\) in \eqref{eq:trace-bound}; the tail bound
controls compactness and the large symmetric rows.  They are
\begin{equation}\label{eq:beta-bounds}
 \begin{gathered}
 \beta_j^{\rm tr\text{-}diff}
 =\frac1{2j(j+1)}+\frac{(-1)^j}4\beta_j^{\rm tail\text{-}sum}
 \quad(j\geq2),\\
 |\beta_j^{\rm tr\text{-}diff}|\leq\frac3{5j^2},\qquad
 0\leq\beta_j^{\rm tail\text{-}sum}\leq\frac3{j^3}
 \quad(j\geq1).
 \end{gathered}
\end{equation}
The first identity follows by subtracting the two trace coefficients
in \eqref{eq:beta}.  For \(j\geq2\), the tail sum has the rational
form
\[
 \beta_j^{\rm tail\text{-}sum}
 =\frac{6(4j^2+4j+3)}{j(j+1)(2j-1)(2j+1)(2j+3)}.
\]
Its required upper bound follows from the exact difference
\[
 \frac3{j^3}-\beta_j^{\rm tail\text{-}sum}
 =\frac{3(12j^3+4j^2-5j-3)}
 {j^3(j+1)(2j-1)(2j+1)(2j+3)}\geq0,
\]
because the numerator polynomial satisfies
\[
 12j^3+4j^2-5j-3
 =12(j-1)^3+40(j-1)^2+39(j-1)+8>0.
\]
For odd \(j\geq3\), the trace difference is positive, since
\[
 \beta_j^{\rm tr\text{-}diff}
 =\frac{4j^3-7j-6}{j(j+1)(2j-1)(2j+1)(2j+3)}>0,
\]
where
\[
 4j^3-7j-6
 =4(j-2)^3+24(j-2)^2+41(j-2)+12.
\]
Its upper bound is already less than \(1/(2j^2)\) by the first
identity in \eqref{eq:beta-bounds}.  For even \(j\geq2\), that
identity makes the trace difference positive, and the upper bound
follows from
\[
 \frac3{5j^2}-\beta_j^{\rm tr\text{-}diff}
 =\frac{4j^4+5j^2-30j-9}
 {5j^2(j+1)(2j-1)(2j+1)(2j+3)}\geq0,
\]
since
\[
 4j^4+5j^2-30j-9
 =4(j-2)^4+32(j-2)^3+101(j-2)^2+118(j-2)+15.
\]
The first mode is checked directly from \eqref{eq:beta}:
\[
 \beta_1^{\rm tr\text{-}diff}=-\frac3{10},\qquad
 \beta_1^{\rm tail\text{-}sum}=\frac15.
\]
The next three trace differences, used in \eqref{eq:trace-bound}, are
\[
 \beta_2^{\rm tr\text{-}diff}=\frac{31}{210},\qquad
 \beta_3^{\rm tr\text{-}diff}=\frac3{140},\qquad
 \beta_4^{\rm tr\text{-}diff}=\frac{157}{4620}.
\]

We now verify the compactness invoked in Proposition~\ref{prop:linear}.
The diagonal and its two neighbors in
\eqref{eq:matrix-upper}--\eqref{eq:matrix-diagonal} are bounded by
\(C/j\), where \(C\) is an absolute constant.  For
\(1\leq i\leq j-2\), \eqref{eq:matrix-dense} and
\eqref{eq:beta-bounds} give
\[
 |K_{ij}|\leq\frac1{4j^2}+\frac{3(2i+1)}{8j^3}
 \leq\frac1{j^2},\qquad
 \sum_{i=1}^{j-2}|K_{ij}|^2\leq\frac1{j^3}.
\]
Together with \eqref{eq:matrix-zero}, these estimates imply
\begin{equation}\label{eq:K-HS}
 \sum_{i,j\geq1}|K_{ij}|^2
 \leq C\sum_{j\geq1}\frac1{j^2}
       +\sum_{j\geq3}\frac1{j^3}<\infty.
\end{equation}
Thus \(\calK\), initially defined on the finite basis span, extends
to a Hilbert--Schmidt and hence compact operator on \(\bbH_\theta\).

\subsection{Symmetric row bounds}

Compactness supplies the range argument in Proposition~\ref{prop:linear}.
Its energy estimate additionally requires the following quantitative
control of the symmetric rows.

\begin{lemma}[Symmetric rows]\label{lem:rows}
The row quantities in \eqref{eq:rows} satisfy
\(\Gamma_1<0\) and \(\Gamma_i<1/16\) for every \(i\geq2\).
\end{lemma}
\begin{proof}
We first isolate the adjacent entries.  The neighboring symmetric
entry has the formula
\begin{equation}\label{eq:adjacent-symmetric}
 A_j:=K_{j-1,j}^{\rm sym}
 =-\frac{2j^3-3j^2-5j-6}
 {2(j+1)(2j-1)(2j+1)(2j+3)}\quad(j\geq2),\qquad A_1=0.
\end{equation}
It is negative for \(j\geq3\), because its numerator equals
\(2(j-3)^3+15(j-3)^2+31(j-3)+6\), while \(A_2=2/105\).
For \(j\geq3\), the identity
\[
 (2j-1)(2j+1)(2j+3)-4(2j^3-3j^2-5j-6)
 =24j^2+18j+21>0
\]
gives \(|A_j|\leq1/(8(j+1))\); the same bound holds at \(j=2\).
Since \(K_{ii}<0\), it remains to estimate the two dense half-rows,
which we denote by
\begin{equation}\label{eq:row-decomposition}
 \begin{gathered}
 \Gamma_i^{\leftarrow}=\frac12\sum_{j=1}^{i-2}|K_{ji}|,
 \qquad
 \Gamma_i^{\rightarrow}=\frac12\sum_{j=i+2}^\infty|K_{ij}|,\\
 \Gamma_i=K_{ii}+|A_i|+|A_{i+1}|
              +\Gamma_i^{\leftarrow}+\Gamma_i^{\rightarrow}.
 \end{gathered}
\end{equation}
Here an empty sum is zero.  The factors \(1/2\) arise because
\eqref{eq:matrix-zero} makes the transposed dense entries vanish.

For \(i\geq3\), \eqref{eq:matrix-dense-trace} and
\eqref{eq:beta-bounds} give
\[
 \begin{aligned}
 \Gamma_i^{\leftarrow}
 &\leq\frac{i-2}{4}|\beta_i^{\rm tr\text{-}diff}|
   +\frac{\beta_i^{\rm tail\text{-}sum}}8
        \sum_{j=1}^{i-2}(j+1)\\
 &\leq\frac3{20i}+\frac3{16i}=\frac{27}{80i},
 \end{aligned}
\]
where \(\sum_{j=1}^{i-2}(j+1)=(i-2)(i+1)/2\leq i^2/2\).
The right half-row satisfies
\[
 \begin{aligned}
 \Gamma_i^{\rightarrow}
 &\leq\frac3{20}\sum_{j=i+2}^\infty\frac1{j^2}
       +\frac{3(i+1)}8\sum_{j=i+2}^\infty\frac1{j^3}\\
 &\leq\frac3{20(i+1)}+\frac3{16(i+1)}
 =\frac{27}{80(i+1)}.
 \end{aligned}
\]
Together with \eqref{eq:row-decomposition}, these bounds give
\begin{equation}\label{eq:large-rows}
 \Gamma_i\leq\frac{27}{80i}+\frac{27}{80(i+1)}
 +\frac1{8(i+1)}+\frac1{8(i+2)}<\frac1{16}\quad(i\geq15).
\end{equation}
The right side is decreasing, and its value at \(i=15\) is
\(6393/108800<1/16\).

For the remaining rows, we use the signs of the dense entries to
sum their absolute values exactly before estimating the result.
Opposite parity in \eqref{eq:matrix-dense} gives a positive entry.
For equal parity, the entry \(K_{ij}\) is negative exactly when
\(j\leq3i\).  Indeed, its sign compares \(2i+1\) with
\[
 \frac{2}{j(j+1)\beta_j^{\rm tail\text{-}sum}}
 =\frac{2j+1}{3}-\frac{2(2j+1)}{(2j+1)^2+2}.
\]
The right side is strictly increasing for \(j\geq2\), is below
\(2i+1\) at \(j=3i\), and is above \(2i+1\) at the next
same-parity index \(j=3i+2\).

In the left half-row, the negative entries therefore occur at
\(j=i-2,i-4,\ldots,i-2d_i^{\rm row}\), where
\(d_i^{\rm row}=\lfloor i/3\rfloor\).  The two arithmetic sums
needed to reverse their signs are
\[
 \sum_{j=1}^{i-2}(-1)^{i+j+1}(2j+1)=1-i+(-1)^i,
 \qquad
 \sum_{k=1}^{d_i^{\rm row}}(2i-4k+1)
 =d_i^{\rm row}(2i-2d_i^{\rm row}-1).
\]
Using \(\sum|K_{ji}|=\sum K_{ji}-2\sum_{K_{ji}<0}K_{ji}\)
therefore gives the exact formula
\begin{equation}\label{eq:left-row}
 \begin{split}
 \Gamma_i^{\leftarrow}
 &=\frac{i-2-2d_i^{\rm row}}{8i(i+1)}
 +\frac{\beta_i^{\rm tail\text{-}sum}}{16}
 \left(1-i+(-1)^i
       +2d_i^{\rm row}(2i-2d_i^{\rm row}-1)\right)
 \quad(i\geq2),
 \end{split}
\end{equation}
with \(\Gamma_1^{\leftarrow}=0\).

For the right half-row, the negative entries are at
\(j=i+2,i+4,\ldots,3i\).  The definition of the tail sum gives
the two telescoping identities
\begin{equation}\label{eq:right-row-telescoping}
 \begin{aligned}
 \sum_{j=i+2}^\infty(-1)^{i+j+1}
        \beta_j^{\rm tail\text{-}sum}
 &=-\frac{\beta_{i+2}^{\rm tail}}{i+2},\\
 \sum_{k=1}^i\beta_{i+2k}^{\rm tail\text{-}sum}
 &=\sum_{j=i+2}^{3i+1}\frac{\beta_j^{\rm tail}}j.
 \end{aligned}
\end{equation}
The first identity sums the signed row; the second collects its
negative terms.  Reversing those signs yields
\begin{equation}\label{eq:right-row-exact}
 \begin{split}
 \Gamma_i^{\rightarrow}
 &=\frac1{8(i+2)}
 -\frac{(2i+1)\beta_{i+2}^{\rm tail}}{16(i+2)}
 +\frac{2i+1}{8}\sum_{j=i+2}^{3i+1}\frac{\beta_j^{\rm tail}}j
 -\frac14\sum_{j=1}^{i}\frac1{(i+2j)(i+2j+1)}.
 \end{split}
\end{equation}
The first remaining sum has the telescoping upper bound
\[
 \begin{aligned}
 \sum_{j=i+2}^{3i+1}\frac{\beta_j^{\rm tail}}j
 &=\sum_{j=i+2}^{3i+1}\frac6{j(4j^2-1)}\\
 &\leq\frac34\sum_{j=i+2}^{3i+1}
 \left(\frac1{(j-1)j}-\frac1{j(j+1)}\right)\\
 &=\frac34\left(\frac1{(i+1)(i+2)}
                   -\frac1{(3i+1)(3i+2)}\right).
 \end{aligned}
\]
For the second, the function \(s\mapsto1/(s(s+1))\) is convex,
and the average of \(i+2,i+4,\ldots,3i\) is \(2i+1\).
Jensen's inequality therefore gives
\[
 \sum_{j=1}^i\frac1{(i+2j)(i+2j+1)}
 \geq\frac{i}{(2i+1)(2i+2)}.
\]
Substitution into \eqref{eq:right-row-exact} proves
\begin{equation}\label{eq:right-row-bound}
 \Gamma_i^{\rightarrow}\leq
 \frac1{8(i+2)}-\frac{(2i+1)\beta_{i+2}^{\rm tail}}{16(i+2)}
 +\frac{3(2i+1)}{32}
 \left(\frac1{(i+1)(i+2)}-\frac1{(3i+1)(3i+2)}\right)
 -\frac{i}{4(2i+1)(2i+2)}.
\end{equation}
Define \(\Gamma_i^{\rm bd}\) as
\(K_{ii}+|A_i|+|A_{i+1}|+\Gamma_i^{\leftarrow}\) plus the
right side of \eqref{eq:right-row-bound}.  By
\eqref{eq:row-decomposition}, it is an upper bound for \(\Gamma_i\).
The following table evaluates this rational expression for the
fourteen remaining indices.  Its first entry is negative, and each
other numerator is less than one sixteenth of its denominator.
\begin{center}
\begin{tabular}{r r@{\qquad}r r}
\toprule
\(i\)&\(\Gamma_i^{\rm bd}\)&\(i\)&\(\Gamma_i^{\rm bd}\)\\
\midrule
1&\(-29/384\)&8&\(752905271/15237331200\)\\
2&\(1787/80640\)&9&\(18127234541/390896150400\)\\
3&\(13297/443520\)&10&\(13028200827/293742310400\)\\
4&\(770179/14414400\)&11&\(5165813479/123992668800\)\\
5&\(103100351/1960358400\)&12&\(63494640661/1597564463040\)\\
6&\(1415487943/26072766720\)&13&\(1984067096231/52743636691200\)\\
7&\(15385927961/299836817280\)&14&\(9573998711/266702374400\)\\
\bottomrule
\end{tabular}
\end{center}
These finite rational comparisons, together with
\eqref{eq:large-rows}, prove the asserted bounds for every row.
\end{proof}

\bibliographystyle{myalpha}
\bibliography{reference}
\end{document}

%% file: macro-3.tex
\usepackage[usenames,dvipsnames]{xcolor}
\usepackage{fancyhdr}
\usepackage{amsmath,amsfonts,amsbsy,amsgen,amscd,amssymb,amsthm}
\usepackage{url}

\usepackage{mathtools}

\mathtoolsset{showonlyrefs}

\usepackage[font=small,margin=0.25in,labelfont={bf},labelsep={colon}]{caption}

\usepackage{subcaption}

\usepackage{tikz}
\usepackage{microtype}
\usepackage{enumitem}

\definecolor{dark-gray}{gray}{0.3}
\definecolor{dkgray}{rgb}{.4,.4,.4}
\definecolor{dkblue}{rgb}{0,0,.5}
\definecolor{medblue}{rgb}{0,0,.75}
\definecolor{rust}{rgb}{0.5,0.1,0.1}

\usepackage[colorlinks=true]{hyperref}

\hypersetup{urlcolor=rust}
\hypersetup{citecolor=dkblue}
\hypersetup{linkcolor=dkblue}

\usepackage{setspace}

\usepackage{graphicx}
\usepackage{booktabs,longtable,tabu} 
\usepackage{multirow} 

\usepackage{float}

\usepackage[full]{textcomp}

\usepackage[scaled=.98,sups,osf]{XCharter}
\usepackage[scaled=1.04,varqu,varl]{inconsolata}
\usepackage[type1]{cabin}
\usepackage[charter,vvarbb,scaled=1.07]{newtxmath}
\usepackage[cal=boondoxo]{mathalfa}
\usepackage[T1]{fontenc}

\usepackage{bm} 

\graphicspath{{figures/}}

\newtheorem{theorem}{Theorem}[section]
\newtheorem{lemma}[theorem]{Lemma}

\newtheorem{proposition}[theorem]{Proposition}

\theoremstyle{definition}

\newtheorem{definition}[theorem]{Definition}

\numberwithin{equation}{section} 
\numberwithin{figure}{section}
\numberwithin{table}{section}

\floatstyle{plaintop}
\newfloat{recipe}{thp}{lor}
\floatname{recipe}{Recipe}
\numberwithin{recipe}{section}

\renewcommand{\phi}{\varphi}

\newcommand{\iunit}{\mathrm{i}}

\newcommand{\Id}{\mathbf{I}}

\newcommand{\diff}{\mathrm{d}}
\newcommand{\idiff}{\, \diff}

\newcommand{\triplenorm}[1]{{\left\vert\kern-0.25ex\left\vert\kern-0.25ex\left\vert #1
    \right\vert\kern-0.25ex\right\vert\kern-0.25ex\right\vert}}

\newcommand{\bbG}{\mathbb{G}}
\newcommand{\bbH}{\mathbb{H}}

\newcommand{\bbR}{\mathbb{R}}

\newcommand{\bbT}{\mathbb{T}}

\newcommand{\bbV}{\mathbb{V}}

\newcommand{\bbX}{\mathbb{X}}
\newcommand{\bbY}{\mathbb{Y}}
\newcommand{\bbZ}{\mathbb{Z}}

\newcommand{\calB}{\mathcal{B}}

\newcommand{\calD}{\mathcal{D}}

\newcommand{\calG}{\mathcal{G}}

\newcommand{\calJ}{\mathcal{J}}
\newcommand{\calK}{\mathcal{K}}
\newcommand{\calL}{\mathcal{L}}

\newcommand{\calN}{\mathcal{N}}

\newcommand{\calP}{\mathcal{P}}

\newcommand{\calR}{\mathcal{R}}
\newcommand{\calS}{\mathcal{S}}
\newcommand{\calT}{\mathcal{T}}

\newcommand{\frA}{\frak A}
\newcommand{\frB}{\frak B}

\newcommand{\frE}{\frak E}
\newcommand{\frF}{\frak F}

\newcommand{\frQ}{\frak Q}

\def\XXint#1#2#3{{%
\setbox0=\hbox{$#1{#2#3}{\int}$}
\vcenter{\hbox{$#2#3$}}\kern-.5\wd0}}

\renewcommand{\leq}{\leqslant}
\renewcommand{\geq}{\geqslant}
\renewcommand{\subset}{\subseteq}

\newcommand{\bmH}{\bm{H}}
